\documentclass{article}
\usepackage{amssymb,amsmath,bm,geometry,graphics,graphicx}
\usepackage{color,url}

\def\no{\noindent}
\def\pmatrix{\left(\begin{array}}
\def\endpmatrix{\end{array}\right)}
\def\udots{\reflectbox{$\ddots$}}  

\def\RR{\mathbb{R}}

\def\D{{\cal D}}

\def\I{{\cal I}}

\def\N{{\cal N}}

\def\P{{\cal P}}

\def\dd{\mathrm{d}}

\newtheorem{theo}{Theorem}
\newtheorem{lem}{Lemma}
\newtheorem{cor}{Corollary}
\newtheorem{rem}{Remark}

\def\proof{\noindent\underline{Proof}\quad}
\def\QED{\mbox{~$\Box{~}$}}

\def\bfb{{\bm{b}}}
\def\bfc{{\bm{c}}}

\def\bfe{{\bm{e}}}

\def\bfuno{{\bm{1}}}

\def\bfgamma{{\bm{\gamma}}}

\def\bfphi{{\bm{\phi}}}

\def\aa{\alpha}

\def\cc{\gamma}
\def\bk{\hat{b}}
\def\ck{\hat{c}}
\def\eps{\varepsilon}

\title{Analysis of Energy and QUadratic Invariant Preserving (EQUIP) methods}

\author{Luigi Brugnano$^a$,~ Gianmarco Gurioli$^a$,~ Felice Iavernaro$^b$\thanks{Corresponding author, felice.iavernaro@uniba.it.} \\
~\\
$^a$ {\small Dipartimento di Matematica e Informatica ``U.\,Dini'', Viale Morgagni 67/a, 50134 Firenze, Italy}\\
$^b$ {\small Dipartimento di Matematica, Via Orabona 4, 70125 Bari, Italy}
}

\begin{document}

\maketitle

\begin{abstract}
In this paper we are concerned with the analysis of a class of  geometric integrators, at first devised in \cite{BIT10-1,BIT12-2}, which can be regarded as an  energy-conserving variant of  Gauss collocation methods.  With these latter they share the property of conserving quadratic first integrals  but, in addition, they also conserve the Hamiltonian function itself. We here reformulate the methods in a more convenient way, and propose a more refined analysis than that given in \cite{BIT12-2} also providing, as a by-product, a practical procedure for their implementation. A thorough comparison with the original Gauss methods is carried out by means of a few numerical tests solving Hamiltonian and Poisson problems.

\medskip
\no{\bf Keywords:} Gauss collocation methods, symplectic methods, energy-conserving methods, line integral methods, Hamiltonian problems, Poisson problems.\\

\smallskip
\no{\bf MSC:} 65P10, 65L05.
\end{abstract}

\section{Introduction}
The study of numerical methods for differential equations  able to retain relevant geometric properties of the solutions, is the core of what nowadays is named {\em geometric integration} (see, e.g., the classical monographs \cite{SSCa94,LeRe04,HLW06} and, more recently, \cite{BI16,BlCa16}). Accordingly, numerical methods able to reproduce specific geometric properties of the continuous dynamical system in the discrete one generated by their application, are often referred to as {\em geometric integrators}. In the framework of ODEs, this approach dates back to the work of G.\,Dahlquist, who devised the linear stability analysis of numerical methods for ODEs  in order to assess their  behavior when solving dissipative problems. In this respect, $A$-stable methods can be considered as geometric integrators, when dealing with such problems. 

For a large class of problems which are of interest in the applications, the geometric properties of the solutions can be taken back to the presence of {\em constants of motion}, which are conserved along the evolution. Such problems are, therefore, named {\em conservative problems}. In more detail, let us consider an ODE-IVP in the form
\begin{equation}\label{fy}
\dot y= f(y), \qquad t>0, \qquad y(0)=y_0\in\Omega\subseteq \RR^m.
\end{equation}
The problem admits $\nu<m$ constants of motion when there exists
\begin{equation}\label{Cy}
C:\Omega \rightarrow \RR^\nu
\end{equation}
such that 
\begin{equation}\label{consC}
y(t) \mbox{~solves~(\ref{fy})} \quad \Rightarrow \quad C(y(t))\equiv C(y_0), \quad \forall t\ge0.
\end{equation}
Hereafter, we shall assume that both $f$ and $C$ are suitably regular.\footnote{E.g., analytical.} From a local point of view, the conservation property (\ref{consC}) is equivalent to require that
\begin{equation}\label{local}\frac{\dd}{\dd t} C(y(t)) = \nabla C(y(t))^\top\dot y(t) = \nabla C(y(t))^\top f(y(t))=0,\end{equation}
i.e., since the initial point of the trajectory is arbitrary,
\begin{equation}\label{local1}
\nabla C(y)^\top f(y)=0, \qquad \forall y\in\Omega.
\end{equation}
An alternative way for requiring the conservation of $C(y(t))$ is through the use of a {\em line integral},
\begin{equation}\label{lineint}
0 = C(y(t)) - C(y(0)) = \int_0^t \nabla C(y(s))^\top\dot y(s)\dd s, \qquad \forall t>0,
\end{equation}
which is clearly  equivalent to (\ref{local}). The reformulation (\ref{lineint}) of the  invariance property is at the basis of the so-called {\em (discrete) Line Integral Methods} (see, e.g., the recent monograph \cite{BI16}),   and here we will exploit this new framework to carry out a theoretical analysis of the methods and to devise a practical implementation strategy.  

The methods  we are going to consider can be regarded as a suitable modification of the basic Gauss-Legendre collocation methods \cite{Bu64}, which we rewrite here by using the well-known $W$-transformation \cite{HW81}. In more details,  introducing the orthonormal Legendre polynomial basis  $\{P_i(x)\}$ on $[0,1]$  ($\Pi_i$ denotes the subspace of polynomials of degree at most $i$)
\begin{equation} \label{Pi} P_i(x)\in\Pi_i, \qquad \int_0^1 P_i(x)P_j(x)\dd x = \delta_{ij}, \qquad \forall i,j=0,1,\dots,
\end{equation}
the matrices
\begin{equation}\label{PO}
\P_s = \pmatrix{ccc} P_0(c_1) &\dots & P_{s-1}(c_1)\\
\vdots & &\vdots\\
P_0(c_s) & \dots &P_{s-1}(c_s) \endpmatrix, \qquad \Omega = \pmatrix{ccc} b_1\\ &\ddots \\ &&b_s\endpmatrix,
\end{equation}
with 
\begin{equation}\label{cibi}
P_s(c_i) = 0, \qquad b_i = \int_0^1 \ell_i(x)\dd x, \qquad \ell_i(x) = \prod_{j\ne i} \frac{x-c_j}{c_i-c_j}, \qquad i=1,\dots,s,
\end{equation}
the Legendre nodes and weights, and matrix
\begin{equation}\label{Xs}
X_s = \pmatrix{cccc} 
\xi_0 & -\xi_1\\
\xi_1 & 0 &\ddots \\
        &\ddots & \ddots & -\xi_{s-1}\\
        &           &\xi_{s-1} &0\endpmatrix,\qquad \xi_i = \frac{1}{2\sqrt{|4i^2-1|}},\quad i=0,\dots,s-1,
        \end{equation}
        the $s$-stage Gauss method is defined by the Butcher tableau
\begin{equation}\label{Butab}
\begin{array}{c|c} \bfc & A:=\P_sX_s\P_s^\top\Omega \\
             \hline\\[-3mm]
             &\bfb^\top\end{array},\qquad \bfb = \pmatrix{c} b_1,\,\dots\,,b_s\endpmatrix^\top, \qquad \bfc = \pmatrix{c} c_1,\,\dots\,,c_s\endpmatrix^\top.
 \end{equation}
             As is well known (see, e.g., \cite{BuBu79}), such methods are algebraically stable, in that
\begin{equation}\label{as}
\Omega A + A^\top \Omega - \bfb\bfb^\top = O,
\end{equation}
and, in particular, they preserve any quadratic invariant of (\ref{fy}). In the framework of geometric integration, condition (\ref{as}) is equivalent to say that such methods are {\em symplectic} \cite{La88,SS88}. The importance of symplectic methods  relies to their application to canonical Hamiltonian systems, i.e. problems of the form
\begin{equation}\label{H}
\dot y = f(y) \equiv J \nabla H(y), \qquad J=\pmatrix{rr} O_{d} &I_{d}\\-I_{d} &O_{d}\endpmatrix = -J^\top = -J^{-1},
\end{equation}
with $m=2d$ and $O_d$ and $I_d$ the zero and identity matrices of dimension $d$, respectively (hereafter $I_r$ will denote the identity matrix of dimension $r$). Moreover,
\begin{equation}\label{Hfun}
H:\Omega\subseteq\RR^m\rightarrow \RR
\end{equation}
is the Hamiltonian function defining the problem, which we assume to be suitably regular. As is well-known, in such a case, the Hamiltonian function is a constant of motion for the dynamical system defined by  (\ref{H})--(\ref{Hfun}).  By virtue of (\ref{as}) a symplectic method does preserve quadratic first integrals including the Hamiltonian function itself in case (\ref{H}) is linear. The capability of conserving quadratic invariants has been recently  recognized as one of the main advantages in using symplectic integrators (see \cite{SS16}).

With this premise, assuming hereafter that only one constant of motion exists for (\ref{fy}) (i.e., $\nu=1$ in (\ref{Cy})), besides possible quadratic invariants,\footnote{We observe that the methods could in principle be generalized to cope with the conservation of more invariants (i.e., $\nu>1$), as sketched in \cite{BIT12-2}. Nevertheless, the arguments become more involved, so that we shall not consider this case, here.}
  the methods we shall study have a Runge-Kutta representation in the form (cf. (\ref{Butab}))
\begin{equation}\label{Butab1}
\begin{array}{c|c} \bfc & A(\aa):=\P_sX_s(\aa)\P_s^\top\Omega \\
             \hline\\[-3mm]
             &\bfb^\top\end{array},
             \end{equation}
where $\P_s$ and $\Omega$ are the same matrices defined in (\ref{PO}), and (see (\ref{Xs})), by setting $\bfe_i\in\RR^s$ the $i$th unit vector,
\begin{equation}\label{Xs1}
X_s(\alpha) := X_s - \alpha W_s \equiv \pmatrix{cccc} 
\xi_0 & \aa-\xi_1\\
\xi_1-\aa & 0 &\ddots \\
        &\ddots & \ddots & -\xi_{s-1}\\
        &           &\xi_{s-1} &0\endpmatrix, \quad W_s = \bfe_2\bfe_1^\top -\bfe_1\bfe_2^\top = -W_s^\top.
\end{equation} 
We observe that, for $\aa=0$, (\ref{Butab1})--(\ref{Xs1}) coincides with the Butcher tableau (\ref{Butab}) of the $s$-stage Gauss collocation method. In the present case,
the scalar parameter $\alpha$ is chosen, at each step of application of the method for numerically solving (\ref{fy}), in order to enforce the conservation of the invariant (\ref{Cy}) (or (\ref{Hfun})) in the numerical solution  (see also \cite{BIT12-2}, for the Hamitonian case).  To this purpose,  exploiting Theorem\,5.1 in \cite[Chapter IV.5]{HW96} and the symmetry\,\footnote{See Theorem \ref{sym} below.} of  the method (\ref{Butab1}), we see that it has order 2 for any fixed $\alpha \not = 0$ and, clearly, order $2s$ when $\alpha=0$. More precisely, denoting by $y_1(\alpha)$ the solution yielded by the application of a single step of method (\ref{Butab1}) to problem (\ref{fy}), we have
\begin{equation}
||y_1(\alpha)-y(h)|| = \alpha O(h^3)+O(h^{2s+1}).
\end{equation}
In Section \ref{equipan} we will show that  an $\alpha^\ast=O(h^{2s-2})$ actually exists such that $C(y_1(\alpha^\ast))=C(y_0)$, which will confer an order $p=2s$ on the resulting method,  as is the case for the corresponding Gauss formula. Since $s\ge 2$, we may confine the study of existence and uniqueness of such an $\alpha^\ast$ inside  a ball of radius $c h^2$. Therefore in the sequel we make the following assumption
\begin{equation}\label{assumption}
|\alpha|\le c h^2, \qquad \mbox{for a given positive constant $c$}.
\end{equation}
The  formulae resulting from the choice $\alpha=\alpha^\ast$ have been named {\em Energy and QUadratic Invariants Preserving} (EQUIP) methods in \cite{BIT12-2} since, beside preserving the Hamiltonian function, they continue to preserve any quadratic invariant just like the Gauss collocation methods they come from. It must be stressed that this  double feature seems to confer great robustness to the methods, as we will see later in the numerical tests.  For sake of completeness, we mention that the idea of modifying Runge-Kutta methods, by introducing a suitable parameter to enforce energy conservation, has been also considered, e.g., in \cite{CHAMR2006,Ko2016}.

An aspect which has not yet been studied neither in  \cite{BIT12-2} nor in the subsequent research concerns the implementation of the methods. Motivated by this issue, in the present work we exploit the characterization (\ref{lineint}) of a constant of motion to show that the methods (\ref{Butab1})  may be indeed interpreted as discrete line integral methods. This will allow us to carry out a more refined theoretical analysis than that given in \cite{BIT12-2} and to devise a practical implementation strategy.

The structure of the paper is as follows: in Section~\ref{equip} we provide full details about the derivation of the parameter $\alpha$; in Section~\ref{equipan} we provide a novel analysis for the corresponding method, besides that given in \cite{BIT12-2}; in Section~\ref{equipdis} we consider the fully discrete method; Section~\ref{fixpoint} is devoted to the solution of the generated discrete problem; Section~\ref{numtest} provides numerical evidence for the theoretical findings; at last, Section~\ref{fine} contains a few conclusions and further directions of investigation.

\section{Energy and QUadratic Invariants Preserving (EQUIP)\\ methods}\label{equip}

In this section, we derive an expression for the parameter $\alpha$ in (\ref{Butab1})--(\ref{Xs1}), chosen in order to enforce the conservation of $C(y)$. Since we are speaking about a one-step method, it will suffice to consider the very first application of the method, in order to obtain the approximation
\begin{equation}\label{y1}
y_1:=y_1(\alpha) \approx y(h), 
\end{equation}
where $h$ is the time step and $y(t)$ is the solution of (\ref{fy}). The derivation of the parameter $\alpha$ will be then derived in the framework of {\em line integral methods}, namely methods which enforce the conservation of the invariants through an approximation of the line integral (\ref{lineint}). Such methods, at first derived for the numerical solution of Hamiltonian problems \cite{IP07,IP08,IT09}, eventually resulted in the class of {\em Hamiltonian Boundary Value Methods (HBVMs)} \cite{BFI14,BIT09,BIS10,BIT10,BIT11,BIT12,BIT12-1,BIT15}, and have been then generalized along several directions \cite{ABI15,BBGCI17,BCMR12,BrFCIa15,BI12,BIT12-3,BS14,YZ2016}, including the methods here studied.\footnote{Despite the common {\em line integral} framework to gain the conservation property, the methods that we shall describe here greatly differ from HBVMs. As matter of fact, the latter methods are not symplectic and are conservative only for Hamiltonian problems, whereas the former ones are defined, at each integration step, by a symplectic map and can be conservative for more general problems.} 

To begin with, we start recalling that, when $\alpha=0$, the method defined by (\ref{Butab1})--(\ref{Xs1}) reduces to the $s$-stage Gauss method, which is defined by the corresponding {\em collocation polynomial} $\sigma\in\Pi_s$ such that
$$\sigma(0) = y_0, \qquad \dot\sigma(c_ih) = f(\sigma(c_ih)),\quad i=1,\dots, s, \qquad y_1\equiv y_1(0) :=\sigma(h).$$
Denoting by $Y_i$, $i=1,\dots,s$ the stages associated with the method,  an equivalent way of defining the collocation polynomial $\sigma(ch)$ is via the interpolation conditions
\begin{equation}\label{path1}
\sigma(0) = y_0, \qquad Y_i = \sigma(c_ih),\quad i=1,\dots, s.
\end{equation}
When $\alpha\ne0$,  the collocation polynomial is replaced with a {\em piecewise continuous polynomial path}, which we continue to denote by $\sigma$,\footnote{For sake of notation simplicity, in the sequel we will  omit to explicitly highlight the dependence of the continuous polynomial path $\sigma$ and some related quantities on the parameter $\alpha$, if this is clear from the context.} 
$$\sigma:[0,h]\cup[h,h+1]\rightarrow \RR^m,$$ such that the interpolation conditions (\ref{path1})
are still fulfilled for the stages $Y_i$ associated with (\ref{Butab1}) and, moreover,
\begin{equation}\label{path2_1}
y_1\equiv y_1(\alpha) := \sigma(h+1), \qquad s.t. \qquad C(\sigma(1+h)) = C(\sigma(0)) ~\Leftrightarrow~ C(y_1)=C(y_0).
\end{equation}
For sake of simplicity, let us denote
\begin{equation}\label{sig12}
\sigma_1(t) \equiv \sigma(t), \quad t\in[0,h],\qquad\quad
\sigma_2(t) \equiv \sigma(t+h), \quad\, t\in[0,1],
\end{equation}
and notice that the interpolation conditions (\ref{path1}) only involve the polynomial $\sigma_1$. The following results hold true.

\begin{theo}\label{s1th} For fixed $\alpha$ and $h$, and with reference to the polynomial basis (\ref{Pi}), let us denote
\begin{equation}\label{gammai}
\cc_j := \sum_{i=1}^s b_i P_j(c_i)f(Y_i), \qquad j=0,\dots,s-1,
\end{equation}
where $Y_i$ are the stages associated with the Runge-Kutta method  (\ref{Butab1}). Moreover, see (\ref{Xs}), we define the vectors
\begin{equation}\label{fii}
\bfphi_i \equiv \pmatrix{c} \phi_{i0}\\ \vdots \\ \phi_{i,s-1}\endpmatrix := X_s^{-1}\bfe_i, \qquad i=1,2.
\end{equation}
Then, the polynomial $\sigma_1\in \Pi_s$ defined as
\begin{equation}\label{sig1}
\sigma_1(ch) = y_0 + h\sum_{j=0}^{s-1} \int_0^c P_j(x)\dd x\,\left[ \cc_j -\aa\left( \phi_{2j} \cc_0 -\phi_{1j}\cc_1\right)\right], \qquad c\in[0,1],
\end{equation}
satisfies the interpolation conditions (\ref{path1}).
\end{theo}
\proof The polynomial $\sigma_1$ defined at (\ref{sig1}) clearly satisfies the first condition in (\ref{path1}). Setting $$Y := \pmatrix{c} Y_1\\ \vdots \\ Y_s\endpmatrix,$$ the stage vector of (\ref{Butab1})--(\ref{Xs1}), and 
\begin{equation}\label{Y}
\sigma_1(\bfc h) := \pmatrix{c}
\sigma_1(c_1h) \\ \vdots \\ \sigma_1(c_sh)\endpmatrix,
\end{equation}
the stage interpolation conditions in (\ref{path1}) read, in compact notation, $Y= \sigma_1(\bfc h)$.
Denoting 
\begin{equation}\label{uno}
\bfuno = (1,\dots,1)^\top\in\RR^s,
\end{equation}
one obtains:
\begin{eqnarray}\nonumber
Y &=& \bfuno\otimes y_0 + h\P_sX_s(\aa)\P_s^\top\Omega \otimes I_m f(Y) 
   ~=~ \bfuno\otimes y_0 + h\P_sX_s\left[I_s-\aa X_s^{-1}W_s\right]\P_s^\top\Omega \otimes I_m f(Y) \\[1mm]  \label{Ydigamma}
   &=& \bfuno\otimes y_0 + h\P_sX_s\left[I_s-\aa(\bfphi_2\bfe_1^\top -\bfphi_1\bfe_2^\top)\right]\P_s^\top\Omega \otimes I_m f(Y).
\end{eqnarray} 
Now, considering the vector form of (\ref{gammai}),
\begin{equation}\label{gammav}
\bfgamma := \pmatrix{c} \gamma_0 \\ \vdots \\ \gamma_s\endpmatrix = \P_s^\top\Omega\otimes I_m f(Y)
\end{equation}
and the property
\begin{equation}\label{Is}
\P_sX_s = \I_s :=\pmatrix{ccc} 
\int_0^{c_1} P_0(x)\dd x& \dots &\int_0^{c_1} P_{s-1}(x)\dd x\\ 
\vdots & &\vdots\\
\int_0^{c_s} P_0(x)\dd x& \dots &\int_0^{c_s} P_{s-1}(x)\dd x\endpmatrix, 
\end{equation}
from (\ref{Ydigamma}) and (\ref{sig1}) we deduce  that
\begin{eqnarray}\label{Ydigamma1}
Y &=& \bfuno\otimes y_0 + h\P_sX_s\left[I_s-\aa(\bfphi_2\bfe_1^\top -\bfphi_1\bfe_2^\top)\right]\otimes I_m\bfgamma
   ~=~ \sigma_1(\bfc h). 
\end{eqnarray}
This concludes the proof. \QED

\medskip
Concerning the polynomial $\sigma_2(t)$ in (\ref{sig12}), it will be used to gain the conservation property. In more details, one has:
\begin{eqnarray}\label{sig20}
\sigma_2(0) &\equiv& \sigma(h) ~=~ \sigma_1(h) ~=~ y_0 + h\left[\gamma_0 -\aa\left( \phi_{20}\gamma_0 - \phi_{10}\gamma_1\right)\right],\\
\sigma_2(1) &\equiv& \sigma(h+1) ~=:~ y_1 ~=~ y_0 + h\sum_{i=1}^s b_i f(Y_i) ~\equiv~ y_0 + h\gamma_0.
\end{eqnarray} 
Consequently, we can choose a linear polynomial:
\begin{equation}\label{sig2}
\sigma_2(t) = y_1 + (t-1)\,\aa h\left( \phi_{20}\gamma_0 - \phi_{10}\gamma_1\right), \qquad t\in[0,1].
\end{equation}
We observe that, from (\ref{sig1}) and (\ref{sig2}), we obtain:\,\footnote{Clearly, from (\ref{sig2})--(\ref{sig112}) one obtains that $\sigma_2(t)\equiv y_1$, when $\aa=0$.}
\begin{equation}\label{sig112}
\dot\sigma_1(ch) = \sum_{j=0}^{s-1} P_j(c)\,\left[ \cc_j -\aa\left( \phi_{2j} \cc_0 -\phi_{1j}\cc_1\right)\right], \qquad \dot\sigma_2(c) \equiv \aa h\left( \phi_{20}\gamma_0 - \phi_{10}\gamma_1\right), \qquad c\in[0,1].~
\end{equation}
Concerning the choice of the parameter $\alpha$ such that the conservation condition $C(y_1)=C(y_0)$ in (\ref{path2_1}) is fulfilled, the following result holds true.

\begin{theo}\label{alfaT1} Let us consider the polynomials $\sigma_1$ and $\sigma_2$ defined at (\ref{sig1}) and (\ref{sig2}), respectively, and the vectors $\gamma_j$ defined at (\ref{gammai}). Morever, we define the vectors
\begin{equation}\label{ros12}
\rho_j(\sigma_1)  = \int_0^1 P_j(c) \nabla C(\sigma_1(ch))\dd c, \quad j=0,\dots,s-1,\qquad
\bar\rho(\sigma_2) = \int_0^1 \nabla C(\sigma_2(c))\dd c.
\end{equation} 
Then $C(y_1)=C(y_0)$, provided that:
\begin{equation}\label{alfa}
\aa = \frac{ \sum_{j=0}^{s-1} \rho_j(\sigma_1)^\top \gamma_j}{\left( \rho_0(\sigma_1)-\bar\rho(\sigma_2)\right)^\top
\left(\phi_{20}\gamma_0-\phi_{10}\gamma_1\right) +\sum_{j=1}^{s-1} \rho_j(\sigma_1)^\top\left(\phi_{2j}\gamma_0-\phi_{1j}\gamma_1\right)}\,.
\end{equation}
\end{theo}
\proof By using the framework provided by line integral methods, from (\ref{sig12}), (\ref{sig1}), (\ref{sig2}),  and (\ref{sig112}), one has:
\begin{eqnarray*}
\lefteqn{ 0 ~=~ C(y_1) - C(y_0) }\\[2mm]
&=& C(\sigma(h+1))-C(\sigma(0)) ~=~ C(\sigma_2(1)) - C(\sigma_2(0)) + C(\sigma_1(h)) - C(\sigma_1(0))\\
&=& \int_0^1 \nabla C(\sigma_2(c))^\top\dot\sigma_2(c)\dd c + \int_0^h \nabla C(\sigma_1(t))^\top\dot\sigma_1(t)\dd t \\ &=& \int_0^1 \nabla C(\sigma_2(c))^\top\dot\sigma_2(c)\dd c + h\int_0^1 \nabla C(\sigma_1(ch))^\top\dot\sigma_1(ch)\dd c\\
&=& \aa h \left[\int_0^1 \nabla C(\sigma_2(c))\dd c\right]^\top \left( \phi_{20}\gamma_0 - \phi_{10}\gamma_1\right) \\
&&~+~ h\int_0^1 \nabla C(\sigma_1(ch))^\top\sum_{j=0}^{s-1} P_j(c)\dd c\,\left[ \cc_j -\aa\left( \phi_{2j} \cc_0 -\phi_{1j}\cc_1\right)\right]\dd c\\
&=& \aa h \bar\rho(\sigma_2)^\top\left( \phi_{20}\gamma_0 - \phi_{10}\gamma_1\right) + h \sum_{j=0}^{s-1}
\rho_j(\sigma_1)^\top\left[ \cc_j -\aa\left( \phi_{2j} \cc_0 -\phi_{1j}\cc_1\right)\right],
\end{eqnarray*}
from which (\ref{alfa}) immediately follows.\QED\medskip

In case of the Hamiltonian problem (\ref{H})--(\ref{Hfun}), the previous result is modified as follows.

\begin{cor}\label{alfaC1} Let us consider the polynomials $\sigma_1$ and $\sigma_2$ defined at (\ref{sig1}) and (\ref{sig2}), respectively, and the vectors $\gamma_j$ defined at (\ref{gammai}). Moreover, with reference to (\ref{H}), let us define the vectors
\begin{equation}\label{gammas12}
\gamma_j(\sigma_1)  = \int_0^1 P_j(c) f(\sigma_1(ch))\dd c, \quad j=0,\dots,s-1,\qquad
\bar\gamma(\sigma_2) = \int_0^1 f(\sigma_2(c))\dd c. 
\end{equation} 
Then $H(y_1)=H(y_0)$, provided that:
\begin{equation}\label{alfaH}
\aa = \frac{ \sum_{j=0}^{s-1} \gamma_j(\sigma_1)^\top J\gamma_j}{\left( \gamma_0(\sigma_1)-\bar\gamma(\sigma_2)\right)^\top J
\left(\phi_{20}\gamma_0-\phi_{10}\gamma_1\right) +\sum_{j=1}^{s-1} \gamma_j(\sigma_1)^\top J\left(\phi_{2j}\gamma_0-\phi_{1j}\gamma_1\right)}\,.
\end{equation}
\end{cor}
\proof The expression
(\ref{alfaH}) follows immediately from (\ref{alfa}) by taking into account (\ref{H}) and setting $C(y)=H(y)$. In fact, from (\ref{ros12}) one obtains: $\gamma_j(\sigma_1) = J\rho_j(\sigma_1)$, $j=0,\dots,s-1$, and $\bar\gamma(\sigma_2) = J\bar\rho(\sigma_2)$.\QED\medskip

The method defined by (\ref{Butab1})--(\ref{Xs1}), with $\aa$ given by (\ref{alfa}) (i.e., (\ref{alfaH}), in the Hamiltonian case (\ref{H})), are named \cite{BIT12-2} $s$-stage EQUIP methods: the name being accounted for by their analysis in the next section.
We conclude this section observing that the discrete line integral tool provides us with a new nonlinear system associated with an EQUIP method, where the stage vector $Y$ is replaced by the vector $\bfgamma$. In fact, inserting (\ref{Ydigamma1}) in (\ref{gammav}) we obtain that the unknown vectors 
$\gamma_j$ are determined by solving the system 
\begin{equation}
\label{bfgamma}
\bfgamma= \P_s^\top\Omega\otimes I_m\, f(\bfuno \otimes y_0+h(\I_s-\alpha \P_sW_s)\otimes I_m\,\bfgamma),
\end{equation}
with $f(\bfuno \otimes y_0+h(\I_s-\alpha \P_sW_s)\otimes I_m\,\bfgamma)$ denoting the block vector containing the application of the function $f$ to each block entry of the argument.  This new representation of an EQUIP method made up of (\ref{bfgamma}), and either (\ref{alfa}) or (\ref{alfaH}), will be  exploited in the sequel to derive a number of theoretical results. Compared with  \cite{BIT12-2}, where only the existence of an $\alpha$ yielding energy conservation was proved without giving details about how $\alpha$ should be searched,  here equation (\ref{alfa})  yields an implicit equation the parameter $\alpha$ must satisfy, and will provide us with a practical strategy for the implementation  of the methods.

\section{Analysis of EQUIP methods}\label{equipan} 
In what follows, we shall obviously assume $s\ge 2$ in (\ref{Butab1})--(\ref{Xs1}).\footnote{In general, $s\ge\nu+1$, in case of $\nu$ constants of motion.}
As already observed, when $\aa=0$, the $s$-stage EQUIP method (\ref{Butab1})--(\ref{Xs1}) reduces to the $s$-stage Gauss-collocation method.
Another relevant property is symmetry.
 
\begin{theo}\label{sym}
For any fixed $\aa\in\RR$ the $s$-stage Runge-Kutta method (\ref{Butab1})--(\ref{Xs1}) is symmetric.
\end{theo} 
\proof Since there are no redundant stages, the symmetry of the method is equivalent to require (see (\ref{uno}))
$$PA(\aa) +A(\aa)P = \bfuno\bfb^\top, \qquad\mbox{and}\qquad P\bfb = \bfb,$$
with 
$$P := \pmatrix{ccc} & &1\\ &\udots\\ 1\endpmatrix=P^\top=P^{-1}\,\in\RR^{s\times s}.$$
The latter condition is trivial. Concerning the former, one has:
\begin{eqnarray*}
PA(\aa)+A(\aa)P&=& PA(0)+A(0)P -\aa\left( P\P_s W_s\P_s^\top\Omega +\P_sW_s\P_s^\top\Omega P\right) \\[2mm]
&=& \bfuno\bfb^\top -\aa\left( P\P_s W_s\P_s^\top\Omega +\P_sW_s\P_s^\top\Omega P\right).
\end{eqnarray*}
Since (see (\ref{PO})) $P\Omega P =\Omega$\, and \,$W_s^\top = -W_s$, we obtain
\begin{eqnarray*}
\left( P\P_s W_s\P_s^\top\Omega +\P_sW_s\P_s^\top\Omega P\right) &=& 
\left( P\P_s W_s\P_s^\top +\P_sW_s\P_s^\top P\right)\Omega\\ &=& 
\left( (P\P_s W_s)\P_s^\top -\P_s(P\P_sW_s)^\top \right)\Omega ~=:~ \Psi.
\end{eqnarray*}
By considering that $P_j(c_i) = (-1)^j P_j(c_{s-i+1})$ and $W_s=\bfe_2\bfe_1^\top -\bfe_1\bfe_2^\top$,
one obtains:
$$ P\P_sW_s ~=~P\P_s\left(\bfe_2\bfe_1^\top -\bfe_1\bfe_2^\top\right) ~=~ - \P_s\left(\bfe_2\bfe_1^\top +\bfe_1\bfe_2^\top\right) ~=:~ -\P_s G,$$ with $G=G^\top$. Consequently, matrix $\P_s G\P_s^\top$ is symmetric and, therefore:
$$\Psi ~=~ -\P_s G\P_s^\top + \P_sG\P_s^\top ~=~O.\QED$$

A further relevant property of the $s$-stage method (\ref{Butab1})--(\ref{Xs1}) is algebraic stability.

\begin{theo}\label{equipT}
For any fixed $\aa\in\RR$ the $s$-stage Runge-Kutta method (\ref{Butab1})--(\ref{Xs1}) satisfies (\ref{as}).
\end{theo}
\proof
In fact, for any fixed $\aa\in\RR$, from (\ref{PO})--(\ref{Xs}) and (\ref{Butab1})--(\ref{Xs1}) one has:
\begin{eqnarray*}
\Omega A(\aa) + A(\aa)^\top \Omega &=& \Omega \P_s(X_s-\aa W_s)\P_s^\top\Omega +\Omega \P_s(X_s-\aa W_s)^\top\P_s^\top\Omega\\
&=& \Omega \P_s \overbrace{(X_s+X_s^\top)}^{=\bfe_1\bfe_1^\top}\P_s^\top\Omega -\aa \Omega \P_s\overbrace{(W_s+W_s^\top)}^{=O}\P_s^\top\Omega ~=~\bfb\bfb^\top +O.\QED
\end{eqnarray*}

As a straightforward consequence, one obtains the following result.

\begin{cor}\label{equipC}
For any fixed $\aa\in\RR$ the $s$-stage Runge-Kutta method (\ref{Butab1})--(\ref{Xs1}) is symplectic and, therefore, conserves all quadratic invariants of (\ref{fy}).
\end{cor}
\proof See \cite{BuBu79,La88,SS88}.\QED

\begin{rem}
The name {\em Energy and QUadratic Invariants Preserving (EQUIP) method} derives from the conservation properties associated with Theorem~\ref{alfaT1} (i.e., Corollary~\ref{alfaC1}, in the Hamiltonian case), and Corollary~\ref{equipC}. 
\end{rem}
 
 We notice that all the properties discussed  in Theorems~\ref{sym} and \ref{equipT}, and in Corollary~\ref{equipC} above,  hold true independently of the choice of the parameter $\alpha$. If we let   equation (\ref{alfa}) (or its variant (\ref{alfaH})) come into play, a primary issue is to ascertain the existence and uniqueness of the solution of the nonlinear system associated with the method, namely (\ref{alfa})--(\ref{bfgamma}) (or (\ref{alfaH})--(\ref{bfgamma})).  For this purpose, we need some additional preliminary results.

\begin{lem}\label{lemma0}
Let $G:\RR\rightarrow V$, vith $V$ a vector space, be a function admitting a Taylor expansion at 0. Then, ~ $$\int_0^1 P_j(c)G(ch)\dd c = O(h^j),\qquad j=0,1,\dots.$$
\end{lem}
\proof See \cite[Lemma~1]{BIT12-1}.\QED

\begin{lem}\label{Odih} With reference to (\ref{gammai}), (\ref{ros12}), and (\ref{gammas12}), one has, under regularity assumptions on $f$ and $C$:
\begin{equation}\label{Odihj}
\gamma_j(\sigma_1) = O(h^j), \qquad \rho_j(\sigma_1)=O(h^j), \qquad \gamma_j = \gamma_j(\sigma_1)+O(h^{2s-j}), \qquad j\ge0,
\end{equation}
and 
\begin{equation}\label{Odih0}
\bar\rho(\sigma_2) = O(1), \qquad \bar\gamma(\sigma_2) = O(1).
\end{equation}
\end{lem}
\proof The thesis easily follows from Lemma~\ref{lemma0}.\QED

\medskip
We need now to discuss the right-hand sides of the expressions at (\ref{alfa}) and (\ref{alfaH}): we shall then study the corresponding numerators and denominators.

\subsection*{Analysis of the numerators in (\ref{alfa}) and (\ref{alfaH})}\label{Nh}
Hereafter, let us denote $\N(h)$ the numerator at the right-hand side of (\ref{alfa}) or (\ref{alfaH}).

\begin{lem}\label{rojgj} Concerning the numerator in (\ref{alfa}), under regularity assumptions on $f$ and $C$, one has
\begin{equation}\label{Odih2s}
\N(h) := \sum_{j=0}^{s-1} \rho_j(\sigma_1)^\top \gamma_j = O(h^{2s}).
\end{equation}
\end{lem}
\proof 
Following the arguments in \cite{BI12}, let us consider the following expansions along the orthonormal basis (\ref{Pi}):
$$\nabla C(\sigma_1(ch)) = \sum_{j\ge0} P_j(c)\rho_j(\sigma_1), \qquad f(\sigma_1(ch)) = \sum_{j\ge0} P_j(c) \gamma_j(\sigma_1), \qquad c\in[0,1],$$
with $\rho_j(\sigma_1)$ and $\cc_j(\sigma_1)$ given by (\ref{ros12}) and (\ref{gammas12}), respectively. Because of (\ref{local1}) and (\ref{Pi}), one has, then:
$$
0 = \int_0^1 \nabla C(\sigma_1(ch))^\top f(\sigma_1(ch))\dd c = \sum_{i,j\ge0} \rho_i(\sigma_1)^\top \gamma_j(\sigma_1) \underbrace{\int_0^1 P_i(c)P_j(c)\dd c}_{=\,\delta_{ij}} = \sum_{j\ge0} \rho_j(\sigma_1)^\top \gamma_j(\sigma_1).
$$
Consequently, by virtue of (\ref{Odihj}), one obtains
$$\sum_{j=0}^{s-1} \rho_j(\sigma_1)^\top \gamma_j(\sigma_1) = -\sum_{j\ge s} \rho_j(\sigma_1)^\top \gamma_j(\sigma_1) = O(h^{2s}).$$
Then,
$$\sum_{j=0}^{s-1} \rho_j(\sigma_1)^\top \gamma_j = \sum_{j=0}^{s-1} \rho_j(\sigma_1)^\top \left[\gamma_j(\sigma_1)+O(h^{2s-j})\right] = O(h^{2s}) + \sum_{j=0}^{s-1} O(h^j)O(h^{2s-j}) = O(h^{2s}). 
\QED$$

For the particular case of Hamiltonian problems, one has the following simpler result.

\begin{lem}\label{gajgj} Concerning the numerator in (\ref{alfaH}), under regularity assumptions on $H$, one has
\begin{equation}\label{Odih2sH}
\N(h) := \sum_{j=0}^{s-1} \gamma_j(\sigma_1)^\top J \gamma_j = O(h^{2s}).
\end{equation}
\end{lem}
\proof In fact, from (\ref{Odihj}), one has, by taking into account that matrix $J$ is skew-symmetric:
$$\sum_{j=0}^{s-1} \gamma_j(\sigma_1)^\top J \gamma_j = \sum_{j=0}^{s-1} \gamma_j(\sigma_1)^\top J \left[\gamma_j(\sigma_1)+O(h^{2s-j})\right] = 0+\sum_{j=0}^{s-1} O(h^j)O(h^{2s-j}) = O(h^{2s}). 
\QED$$

\subsection*{Analysis of the denominators in (\ref{alfa}) and (\ref{alfaH})}\label{Dh}

Before studying the denominators at the right-hand sides of (\ref{alfa}) and (\ref{alfaH}), which we shall hereafter denote by $\D(h)$, we need to state the following result concerning  the leading entries of the vectors defined at (\ref{fii}).

\begin{lem}\label{fiiL}
With reference to (\ref{Xs}) and (\ref{fii}), one has (discarding $\phi_{ij}$, when $j\ge s$):
\begin{eqnarray}\label{fi1}
\pmatrix{cccc} \phi_{10}, &\phi_{11}, & \phi_{12},& \phi_{13}\endpmatrix &=& \left\{\begin{array}{cl}
(\,\xi_0^{-1},\, 0,\, \xi_1(\xi_0\xi_2)^{-1},\, 0\,), &~ s \quad\mbox{odd},\\[3mm]
(\,0,\, -\xi_1^{-1},\, 0,\, \ne0\,), &~ s \quad\mbox{even},\end{array}\right.\\[3mm] \label{fi2}
\pmatrix{cccc} \phi_{20}, &\phi_{21}, & \phi_{22},& \phi_{23}\endpmatrix &=& \left\{\begin{array}{cl}
(\,0,\,0,\, -\xi_2^{-1},\,0\,), &~ s \quad\mbox{odd},\\[3mm]
(\,\xi_1^{-1},\,\xi_0\xi_1^{-2},\,0,\,\ne0\,), &~ s \quad\mbox{even}.\end{array}\right.
\end{eqnarray}
\end{lem}

Moreover, it is quite straightforward, though lengthy, to prove the following expansions, which hold under the assumption (\ref{assumption}).

\begin{lem}\label{expan} Let us consider the vectors defined at (\ref{gammai}), (\ref{ros12}), and (\ref{gammas12}),  and the parameter $\aa$ defined at either  (\ref{alfa}) or (\ref{alfaH}). Moreover, let us set
\begin{eqnarray}\nonumber
&f_0:=f(y_0), \qquad f'_0:=f'(y_0), \qquad f''_0:=f''(y_0),& \\ 
\label{notafg}\\[-3mm] \nonumber
&g(y):=\nabla C(y) \qquad g_0:=g(y_0), \qquad g'_0:=g'(y_0), \qquad g''_0:=g''(y_0).&
\end{eqnarray}
Then, under regularity assumptions on the involved functions, and assuming that (\ref{assumption}) is satisfied, the following expansions hold true
for the following coefficients defined at (\ref{gammai}), (\ref{ros12}), and (\ref{gammas12}):\,\footnote{They are obtained by expanding $f(y)$ and $\nabla C(y)$ at $y=y_0$.}
\begin{eqnarray} \nonumber
 \cc_0(\sigma_1),\cc_0 &=& f_0+\xi_0hf'_0f_0+h^2[(\xi_0^2-\xi_1^2)(f'_0)^2f_0+\xi_0^{-1}\xi_1^2f''_0(f_0,f_0)]+O(h^3), \\[1mm]  \nonumber
\cc_1(\sigma_1),\cc_1 &=&  \xi_1hf'_0f_0+\xi_0\xi_1h^2[(f'_0)^2f_0+f''_0(f_0,f_0)] +O(h^3),\\[1mm]  \label{expgam}
\cc_2(\sigma_1),\cc_2 &=&  \xi_1\xi_2h^2[(f'_0)^2f_0+f''_0(f_0,f_0)]+O(h^3),\\[1mm]   \nonumber
\bar\cc(\sigma_2) &=& f_0+hf'_0f_0+\xi_0h^2[(f'_0)^2f_0+f''_0(f_0,f_0)]+O(h^3),  
\end{eqnarray}
and 
\begin{eqnarray}\nonumber
\rho_0(\sigma_1) &=&    g_0+\xi_0hg'_0f_0+h^2[(\xi_0^2-\xi_1^2)g'_0f'_0f_0+\xi_0^{-1}\xi_1^2g''_0(f_0,f_0)]+O(h^3),\\[1mm] \nonumber
\rho_1(\sigma_1) &=& \xi_1hg'_0f_0+\xi_0\xi_1h^2[g'_0f'_0f_0+g''_0(f_0,f_0)]+O(h^3),\\[1mm] \label{expro}
\rho_2(\sigma_1) &=&  \xi_1\xi_2h^2[g'_0f'_0f_0+g''_0(f_0,f_0)] +O(h^3),\\[1mm] \nonumber
\bar\rho(\sigma_2) &=&  g_0+hg'_0f_0+\xi_0h^2[g'_0f'_0f_0+g''_0(f_0,f_0)]+O(h^3).
\end{eqnarray}
\end{lem}

The next results then follow.

\begin{lem}\label{Odih2} Let us consider the denominator at the right-hand side in (\ref{alfa}).
Then,  under regularity assumptions on $f$ and assuming (\ref{assumption}) be satisfied, one has:
\begin{equation}\label{D1D2}
\D(h)  = \xi_1h^2\left[f_0^\top((g'_0)^\top f'_0-(g'_0f'_0)^\top)f_0-(g''_0(f_0,f_0))^\top f_0\right]+O(h^3).
\end{equation}
\end{lem}
\proof We set \,$\D(h)=D_1(h)+D_2(h)$,\, with
\begin{eqnarray*} \label{D1}
D_1(h) &:=&\left( \rho_0(\sigma_1)-\bar\rho(\sigma_2)\right)^\top
\left(\phi_{20}\gamma_0-\phi_{10}\gamma_1\right)+\rho_1(\sigma_1)^\top\left(\phi_{21}\gamma_0-\phi_{11}\gamma_1\right), \\ \label{D2}
D_2(h) &:=& \sum_{j=2}^{s-1} \rho_j(\sigma_1)^\top\left(\phi_{2j}\gamma_0-\phi_{1j}\gamma_1\right).
\end{eqnarray*}
Then,  by taking into account (\ref{Odihj}) and (\ref{fi1})--(\ref{expro}), we discuss separately the case $s$ odd and $s$ even.
\begin{description}
\item{\underline{$s$ odd:}} in such a case, one has, by taking into account that $\phi_{22}\xi_1\xi_2=-\xi_1$,
\begin{eqnarray*}
D_1(h)&=&-\xi_0^{-1}(\rho_0(\sigma_1)-\bar{\rho}(\sigma_2))^{\top}\gamma_1~=~-\xi_0^{-1}[-\xi_0hg'_0f_0+O(h^2)]^{\top}(\xi_1hf'_0f_0+O(h^2))\\[1mm]
&=&\xi_1h^2(g'_0f_0)^{\top}(f'_0f_0)~+~O(h^3),\\[2mm]
D_2(h) &=& \phi_{22}\rho_2(\sigma_1)^\top\gamma_0 +O(h^3) ~=~ \phi_{22}\xi_1\xi_2 h^2[g'_0f'_0f_0+g''_0(f_0,f_0)]^{\top}f_0~+~O(h^3)\\[1mm]
&=&-\xi_1h^2[(g'_0f'_0f_0)^{\top}f_0+(g''_0(f_0,f_0))^{\top}f_0]~+~O(h^3).
\end{eqnarray*}
Consequently, (\ref{D1D2}) follows.

\item{\underline{$s$ even:}}  in such a case, one has, by taking into account that ~$(\xi_0^{-1}\xi_1-\xi_0\xi_1^{-1})=-\xi_0^{-2}\xi_1$~ and ~$(\xi_0^2\xi_1^{-1}-\xi_0^{-2}\xi_1)=-\xi_1$,
\begin{eqnarray*}
D_1(h) &=& \xi_1^{-1}(\rho_0(\sigma_1)-\bar{\rho}(\sigma_2))^{\top}\gamma_0+\rho_1(\sigma_1)^{\top}(\xi_0\xi_1^{-2}\gamma_0+\xi_1^{-1}\gamma_1)\\[1mm]
&=&\xi_1^{-1}[-\xi_0hg'_0f_0-\xi_0^{-2}\xi_1^2h^2g'_0f'_0f_0+(\xi_0^{-1}\xi_1^2-\xi_0)h^2g''_0(f_0,f_0)]^{\top}(f_0+\xi_0hf'_0f_0)~+~\\[1mm]
&&\xi_0\xi_1^{-2}[\xi_1hg'_0f_0+\xi_0\xi_1h^2[g'_0f'_0f_0+g''_0(f_0,f_0)]]^{\top}(f_0+\xi_0hf'_0f_0)~+~\\[1mm]
&&\xi_1^{-1}(\xi_1hg'_0f_0)^{\top}(\xi_1hf'_0f_0)~+~O(h^3)\\[1mm]
&=&-\xi_1h^2(g'_0f'_0f_0)^{\top}f_0-\xi_1h^2(g''_0(f_0,f_0))^{\top}f_0+
\xi_1h^2(g'_0f_0)^{\top}(f'_0f_0)+O(h^3), \\ [2mm]D_2(h) &=& \phi_{23}\rho_3(\sigma_1)^\top\gamma_0 +O(h^4) ~=~ O(h^3),
\end{eqnarray*}
from which (\ref{D1D2}) follows.\QED
\end{description}

In the Hamiltonian case, by considering (\ref{H})--(\ref{Hfun})~ (i.e., $J g(y) \equiv f(y)$),  and the expansions (\ref{expgam}), one proves the following result in a similar way.

\begin{lem}\label{Odih21} Let us consider  the denominator at the right-hand side in (\ref{alfaH}).
Then, under regularity assumptions on $f$ and assuming (\ref{assumption}) be satisfied, one has:
\begin{equation}\label{D1D2H}
\D(h) = -\xi_1h^2 \left[f''_0(f_0,f_0)\right]^\top J f_0 \,+\, O(h^3).  
\end{equation}
\end{lem}

\section*{Existence and uniqueness}\label{EU}

As a simple consequence of the previous Lemmas~\ref{rojgj}--\ref{Odih21}, one obtains the following result, which provides a different, and more general, proof than what proved in \cite[Theorem~3.5]{BIT12-2} in the Hamiltonian case.

\begin{theo}\label{Odih2s2} Assume that $\aa$ satisfies equation  (\ref{alfa}) and
\begin{equation}\label{non0}
\left[f_0^\top((g'_0)^\top f'_0-(g'_0f'_0)^\top)f_0-(g''_0(f_0,f_0))^\top f_0\right]\ne0. 
\end{equation}
Then, \begin{equation}\label{alfaOdih} \aa = O(h^{2(s-1)}).\end{equation}
The same result holds true in the Hamiltonian case, when $\aa$ is given by the right-hand side in  (\ref{alfaH}), and
\begin{equation}\label{non0H}
\left[f''_0(f_0,f_0)\right]^\top J f_0\ne0.
\end{equation}
\end{theo}
\proof Depending on the two cases, one obtains $\aa=\N(h)/\D(h)$, with $\N(h)$ and $\D(h)$ defined according to
(\ref{Odih2s}) and (\ref{D1D2}), when $\aa$ is given by (\ref{alfa}), or (\ref{Odih2sH}) and (\ref{D1D2H}), 
when $\aa$ is given by (\ref{alfaH}). In both cases, (\ref{alfaOdih}) easily follows, due to the fact that $\N(h)=O(h^{2s})$ and $\D(h)=O(h^2)$.\QED

\medskip
\begin{rem} We observe that when (\ref{non0}) is not true, then the leading term in (\ref{D1D2}) vanishes. Practically, this event can be recognized by the slow down or failure of the convergence of the nonlinear iteration later described in Section~\ref{fixpoint} which, in turn, relies on the result of Theorem~\ref{exists_unique} below. In such a case, one may set $\aa=0$, thus performing a single step of the basic $s$-stage Gauss method. The same argument applies, in the Hamiltonian case, when (\ref{non0H}) is not satisfied, 
so that the leading term in (\ref{D1D2H}) vanishes. This strategy will allow to retain the order $2s$ of the underlying $s$-stage Gauss method, as we shall see later
in the analysis of the method.\end{rem}

We are now in the right position to prove the existence and uniqueness result. The nonlinear system to be solved at each step of the integration procedure is  (see (\ref{alfa})--(\ref{alfaH}) and (\ref{bfgamma}))\,\footnote{We now emphasize, in the first equation in (\ref{nonlsys}), the dependence of the numerator and denominator defining $\aa$ from the given arguments, rather than from the stepsize $h$.}
\begin{equation}
\label{nonlsys}
\left\{ 
\begin{array}{l}\displaystyle \alpha=\frac{N(\aa,\bfgamma)}{D(\aa,\bfgamma)} := \frac{\N(h)}{\D(h)},\\[.4cm]
\bfgamma=\Psi(\aa,\bfgamma):= \P_s^\top\Omega\otimes I_m\, f(\bfuno \otimes y_0+h(\I_s-\alpha \P_sW_s)\otimes I_m\,\bfgamma).
\end{array}
\right.
\end{equation}  
Setting $\alpha=0$ in the first equation of (\ref{nonlsys}) leads back to the nonlinear system associated with the Gauss method of order $2s$ that we write as 
\begin{equation}
\label{nonlsysGauss}
\bfgamma = \Psi_0(\bfgamma) := \Psi(0,\bfgamma). 
\end{equation}
Assuming $f:\RR^m \rightarrow \RR^m$ Lipschitz continuous with constant $L>0$, we see that $\Psi_0(\bfgamma)$ is Lipschitz continuous in $\RR^{sm}$ with constant\,\footnote{Hereafter, $\|\cdot\|$ denotes the infinity norm.} 
\begin{equation}\label{K0}
K_0:=hL\|\P_s^\top\Omega\|\,\|\I_s\|=O(h),
\end{equation}
therefore $\Psi_0(\bfgamma)$ is a contraction on $\RR^{sm}$ for $h$ sufficiently small, and the Banach fixed-point theorem assures the existence of a unique solution $\bar \bfgamma$ of (\ref{nonlsysGauss})  which can be determined as the limit of the sequence $\bfgamma^{\ell+1}=\Psi_0(\bfgamma^\ell)$, $\ell=0,1,\dots,$ starting from any $\bfgamma_0\in \RR^{sm}$.

\begin{theo}\label{exists_unique}  
Under regularity assumptions on the functions $f$,  and assuming either (\ref{non0}) or (\ref{non0H}) holds true, so that  $D(\aa,\bfgamma)$ is bounded away from zero in a neighbourhood of $(0,\bar \bfgamma)$, system (\ref{nonlsys}) admits a unique solution $(\alpha^\ast,\bfgamma^\ast)$  which can be found as the limit of the sequence 
\begin{equation}
\label{iteration1}
\left\{ 
\begin{array}{l}
\displaystyle \alpha^{\ell+1}=\frac{N(\aa^\ell,\bfgamma^\ell)}{D(\aa^\ell,\bfgamma^\ell)}, \\[.4cm]
\bfgamma^{\ell+1}=\Psi(\aa^\ell,\bfgamma^\ell),  \qquad \ell=0,1,\dots,
\end{array}
\right.
\end{equation}
starting at $$(\aa^0,\bfgamma^0)=(0,\bar \bfgamma).$$ Moreover the solution $(\alpha^\ast,\bfgamma^\ast)$ satisfies 
\begin{equation}
\label{errgamma}
 |\alpha^\ast| =O(h^{2s-2}), \qquad \|\bfgamma^\ast - \bar \bfgamma\|=O(h^{2s-1}).
\end{equation}
\end{theo}
\proof In accord with the contraction mapping theorem, we show that constants $h_0>0$ and $r>0$ exist such that, for $h\in[0,h_0)$, the iteration function, 
\begin{equation}
\label{iterF}
F(\aa,\bfgamma):=\pmatrix{c}
N(\aa,\bfgamma)/D(\aa,\bfgamma)\\[1mm] \Psi(\aa,\bfgamma)\endpmatrix :\RR^{sm+1}\rightarrow \RR^{sm+1},
\end{equation}
defined at (\ref{iteration1}) satisfies the two properties:
\begin{itemize}
\item[(a)] $F(\aa,\bfgamma)$ is  a contraction in the closed ball $B((0,\bar \bfgamma), r)$ of center $(0,\bar \bfgamma)$ and radius $r>0$;
\item[(b)] the sequence $(\aa^\ell,\bfgamma^\ell)\in B((0,\bar \bfgamma), r)$,\, for all\, $\ell=0,1,\dots$.
\end{itemize}
Expanding $\Psi(\aa,\bfgamma)$ with respect to the variable $\alpha$ in a neighbourhood of zero yields
$$
\Psi(\aa,\bfgamma)=\Psi_0(\gamma)+R(\aa,\gamma),
$$
with
$$
R(\aa,\gamma):=-h\alpha\, \P_s^\top\Omega\otimes I_m \, f'(\bfuno \otimes y_0+h\I_s\otimes I_m\bfgamma)\,\P_sW_s\otimes I_m\,\bfgamma ~+~ O((\alpha h)^2),
$$ 
where $f'$ stands for the Jacobian matrix of $f$.
The function $R(\aa,\bfgamma)$ satisfies $R(0,\bfgamma)=0$ and, since by virtue of Theorem~\ref{Odih2s2} we may assume $|\alpha| \le c_1 h^{2s-2}$, we also get 
\begin{equation}
\label{boundR}
\|R(\aa,\gamma)\|\le c_2 h^{2s-1}, 
\end{equation} 
with $c_1$ and $c_2$  positive constants independent of $h$, $\alpha$, and $\bfgamma$, in a closed ball $B((0,\bar \bfgamma), r)$ of center $(0,\bar \bfgamma)$ and given radius $r>0$. Consequently, in a neighbourhood of $(0,\bar \bfgamma)$, system  (\ref{nonlsys}) may be regarded as a perturbation of system (\ref{nonlsysGauss}), the perturbing term being the function $R(\aa,\bfgamma)$. 

From (\ref{boundR}), it follows that $R(\aa,\bfgamma)$ is Lipschitz continuous with constant  $K_1=O(h^{2s-1})$ and hence $\Psi(\aa,\bfgamma)$ is Lipschitz continuous with constant $\mu \le K_0+K_1=O(h)$,  where $K_0$ is defined in (\ref{K0}). This property, together with the Lipschitz continuity of $N(\aa,\bfgamma)/D(\aa,\bfgamma)$ with constant $O(h^{2s-2})$, assures that $h_0>0$ exists such that, for $h\in[0,h_0)$, the iteration function $F(\aa,\bfgamma)$ defined at (\ref{iterF}) is indeed a contraction in $B((0,\bar \bfgamma), r)$, with constant $\mu$.  

To show that system  (\ref{nonlsys}) actually admits a solution $\bfgamma^\ast$ which is $O(h^{2s-1})$-close to $\bar \bfgamma$, we set $$r:=|N(0,\bar \bfgamma)/D(0,\bar \bfgamma)|/(1-\mu)=O(h^{2s-2}),$$ and prove that the sequence $(\aa^\ell,\bfgamma^\ell)$ defined at (\ref{iteration1}) is entirely contained in the closed ball $B((0,\bar \bfgamma), r)$. We proceed by induction. For $\ell=1$ we get
$$
\|\bfgamma^1-\bar \bfgamma\| = \|\Psi(0,\bfgamma^0)-\bar \bfgamma\| = \|\Psi(0,\bar \bfgamma)-\bar \bfgamma\| =  \|\Psi_0(\bar \bfgamma)-\bar \bfgamma\| = 0,
$$
and
$$
|\alpha^1-\alpha^0| = |\alpha^1|= \left|\frac{N(0,\bfgamma^0)}{D(0,\bfgamma^0)}\right| = \left|\frac{N(0,\bar \bfgamma)}{D(0,\bar \bfgamma)}\right| = (1-\mu) r. 
$$
Assuming now that $(\aa^\ell,\bfgamma^\ell) \in B((0,\bar \bfgamma),r)$, we obtain
$$
\|\bfgamma^{\ell+1}-\bar \bfgamma\| =  \|\Psi(\aa^\ell,\bfgamma^\ell)- \Psi(0,\bar \bfgamma)\| \le \mu \max\{\|\bfgamma^\ell-\bar \bfgamma\|, |\alpha^\ell|\} \le \,\mu r \,<r,
$$
and
$$
\begin{array}{rcl}
|\alpha^{\ell+1}-\alpha^0| &\le&  \displaystyle|\alpha^{\ell+1}-\alpha^1|+|\alpha^{1}-\alpha^0| = \left|\frac{N(\aa^\ell,\bfgamma^\ell)}{D(\aa^\ell,\bfgamma^\ell)}- \frac{N(\aa^0\bfgamma^0)}{D(\aa^0,\bfgamma^0)} \right| + |\alpha^{1}| \\[.3cm]
&\le& \displaystyle \mu \max\{\|\bfgamma^\ell-\bar \bfgamma\|, |\alpha^\ell|\} + |\alpha^{1}| \le \mu r+  (1-\mu) r = r. 
\end{array}
$$
Since $(\aa^\ast,\bfgamma^\ast) \in B((0,\bar \bfgamma),r)$, then clearly $|\aa^\ast| \le r =O(h^{2s-2})$. Moreover,
subtracting $\bar \bfgamma=\Psi(0,\bar \bfgamma)$ from $\bfgamma^\ast=\Psi(\aa^\ast,\bfgamma^\ast)$ we finally obtain
$$
\|\bfgamma^\ast-\bar \bfgamma\|=\|\Psi(\alpha^\ast,\bfgamma^\ast)-\Psi(0,\bar \bfgamma)\|\le  \mu \max\{\|\bfgamma^\ast-\bar \bfgamma\|, |\alpha^\ast|\}=O(h^{2s-1}). \QED
$$


\section*{Order of convergence}

Next, we show that method (\ref{Butab1})--(\ref{Xs1}), with $\aa$ given by (\ref{alfa}) (or (\ref{alfaH}), in case the problem is Hamiltonian), retains the order $2s$ of the underlying $s$-stage Gauss method (i.e., the method corresponding to $\aa=0$). 

\begin{theo}\label{order} Let $y_1\approx y(h)$ be the approximation provided by the $s$-stage EQUIP method (\ref{Butab1})--(\ref{Xs1}), with $\aa$ given by (\ref{alfa}) (respectively, (\ref{alfaH}) when the problem is Hamiltonian). Then,  assuming that the hypotheses of Theorem~\ref{exists_unique} hold true,
\begin{equation}\label{ord2s}
y_1-y(h) = O(h^{2s+1}), 
\end{equation}
i.e., the method has order $2s$.
\end{theo}
\proof The proof will follow the novel approach defined in \cite{BIT12-1}. For this purpose, let us denote by
$y(t,\theta,\xi)$ the solution of the problem $$\dot y=f(y), \qquad t>\theta, \qquad y(\theta)=\xi.$$ Moreover, let $\Phi(t,\theta)$ be the fundamental matrix solution of the corresponding variational problem 
$$\dot \Phi(t,\theta) = f'(y(t,\theta,\xi)) \Phi(t,\theta), \qquad t>\theta, \qquad \Phi(\theta,\theta)=I_m.$$
We also recall that, according to Lemma~\ref{lemma0},
\begin{equation}\label{Gaj}
\Gamma_j(h) := \int_0^1 P_j(c) \Phi(h,ch) \dd c ~=~ O(h^j), \qquad j\ge0.
\end{equation}
By taking into account (\ref{path1})--(\ref{sig2}), let us then define the  polynomial
\begin{equation}\label{uu1}
u(ch) = y_0 + h\sum_{j=0}^{s-1} \int_0^c P_j(x)\dd x\,\gamma_j \quad\left[\quad \Rightarrow\quad 
\dot u(ch) = \sum_{j=0}^{s-1} P_j(c)\gamma_j\right], \qquad c\in[0,1],
\end{equation}
such that $u(0)=y_0$ and $u(h)=y_1$. Consqeuently, 
\begin{eqnarray}\nonumber
y_1-y(h)&=& y(h,h,u(h))-y(h,0,u(0)) = \int_0^h \frac{\dd}{\dd t} y(h,t,u(t))\dd t \\ \nonumber
&=& \int_0^h \left[\left.\frac{\partial}{\partial \theta} y(h,\theta,u(t))\right|_{\theta=t} + 
\left.\frac{\partial}{\partial \xi} y(h,t,\xi)\right|_{\xi=u(t)}\dot u(t)\right]\dd t\\ 
&=& \int_0^h \Phi(h,t)\left[ -f(u(t)) + \dot u(t)\right]\dd t = h\int_0^1 \Phi(h,ch)\left[ -f(u(ch)) + \dot u(ch)\right]\dd c \,=: \label{(a)}
\end{eqnarray}
By considering that
$$
f(u(ch)) = \sum_{j\ge 0} P_j(c) \gamma_j(u), \qquad \gamma_j(u) = \int_0^1 P_j(\tau)f(u(\tau h))\dd\tau ~=~ O(h^j),
$$
and taking into account  (\ref{Odihj}), (\ref{Gaj}), and (\ref{uu1}), one  obtains:
\begin{eqnarray*}
(\ref{(a)}) &=& h\int_0^1 \Phi(h,ch)\left[ \sum_{j=0}^{s-1} P_j(c)(\gamma_j(\sigma_1)-\gamma_j(u)) - \sum_{j\ge s}P_j(c)\gamma_j(u) + \sum_{j=0}^{s-1} P_j(c) (\gamma_j-\gamma_j(\sigma_1))\right]\dd c\\
&=& h\left[ \sum_{j=0}^{s-1} \Gamma_j(h)(\gamma_j(\sigma_1)-\gamma_j(u)) - \sum_{j\ge s} \Gamma_j(h)\gamma_j(u) + \sum_{j=0}^{s-1}\Gamma_j(h) (\gamma_j-\gamma_j(\sigma_1))\right].
\end{eqnarray*}
Due to the fact that $\Gamma_j(h) = O(h^j)$, $\gamma_j(u) = O(h^j)$, $\gamma_j-\gamma_j(\sigma_1) = O(h^{2s-j})$, one then concludes that
\begin{equation}\label{erro1}
y_1-y(h) ~=~ h\sum_{j=0}^{s-1} \Gamma_j(h)(\gamma_j(\sigma_1)-\gamma_j(u)) \,+\, O(h^{2s+1}).
\end{equation}
Consequently, (\ref{ord2s}) follows provided that 
\begin{equation}\label{somma}
\sum_{j=0}^{s-1} \Gamma_j(h)(\gamma_j(\sigma_1)-\gamma_j(u)) ~=~ O(h^{2s}).
\end{equation}
By taking into account (\ref{sig1}), (\ref{uu1}), and (\ref{alfaOdih}), one has
\begin{equation}\label{s1menou}
\sigma_1(ch)-u(ch) = -h\aa\sum_{i=0}^{s-1} \int_0^c P_i(x)\dd x\left( \phi_{2i} \cc_0 -\phi_{1i}\cc_1\right) = O(h^{2s-1}),
\end{equation}
so that:
\begin{eqnarray*}
\gamma_j(\sigma_1)-\gamma_j(u) &=& \int_0^1 P_j(c)\left[ f(\sigma_1(ch))-f(u(ch))\right]\dd c \\
&=& \int_0^1 P_j(c)f'(u(ch))\left( \sigma_1(ch)-u(ch) \right)\dd c \,+\, O(h^{4s-2})\\
&=& -h\aa f'(y_0) \int_0^1 P_j(c)\sum_{i=0}^{s-1} \int_0^c P_i(x)\dd x\left( \phi_{2i} \cc_0 -\phi_{1i}\cc_1\right)\dd c ~+~O(h^{2s}) \\
&=& -h\aa f'(y_0)\int_0^1 g_j(c) \dd c ~+~ O(h^{2s}),
\end{eqnarray*}
having set
$$g_j(c) := P_j(c)\sum_{i=0}^{s-1} \int_0^c P_i(x)\dd x\left( \phi_{2i} \cc_0 -\phi_{1i}\cc_1\right).$$
Since $g_j$ is a polynomial of degree $s+j$, for  $j=0,\dots,s-1$, one has that
$$\int_0^1 g_j(c)\dd c = \sum_{\ell = 1}^s b_\ell g_j(c_\ell).$$ Consequently,  by setting
$$\Gamma(h) := \pmatrix{ccc} \Gamma_0(h), &\dots\,,&\Gamma_{s-1}(h)\endpmatrix,$$
and taking into account (\ref{PO}), (\ref{Is}), (\ref{fii}), (\ref{gammav}), one has:
\begin{eqnarray}\nonumber
\lefteqn{
\sum_{j=0}^{s-1} \Gamma_j(h)(\gamma_j(\sigma_1)-\gamma_j(u))}\\ \nonumber
&=&-h\aa \sum_{j=0}^{s-1} \Gamma_j(h) f'(y_0) \sum_{\ell=1}^s b_\ell 
P_j(c_\ell)\sum_{i=0}^{s-1} \int_0^{c_\ell} P_i(x)\dd x\left( \phi_{2i} \cc_0 -\phi_{1i}\cc_1\right)
\,+\,O(h^{2s})\\[1mm] \nonumber
&=& -h\aa \,\Gamma(h)\, (I_s\otimes f'(y_0))\, \left(\P_s^\top\Omega \I_s\otimes I_m\right)\, \left[X_s^{-1}\bfe_2\otimes \cc_0-X_s^{-1}\bfe_1\otimes \cc_1\right]~+~ O(h^{2s})\\[1mm]
&=& -h\aa \,\Gamma(h)\, (I_s\otimes f'(y_0))\, \left(\P_s^\top\Omega \P_s\otimes I_m\right)\, \left[(\bfe_2\bfe_1^\top-\bfe_1\bfe_2^\top)\otimes I_m \right]\, \bfgamma ~+~ O(h^{2s}) \,\quad =:\label{(b)}
\end{eqnarray}
By considering that $\P_s^\top\Omega \P_s=I_s$, one then obtains:
\begin{equation}\label{erro2}
(\ref{(b)}) ~=~ h\aa\left[ \,\Gamma_0(h)f'(y_0) \gamma_1 \,-\,\Gamma_1(h) f'(y_0) \gamma_0 \,\right] ~+~ O(h^{2s}).
\end{equation}
(\ref{somma}) then follows by taking into account that $\aa=O(h^{2s-2})$, because of Theorem~\ref{exists_unique}, and, moreover,  
$$\gamma_1=O(h)=\Gamma_1(h), \qquad \gamma_0=O(1)=\Gamma_0(h).\QED$$

\section{Discretization of the involved integrals}\label{equipdis}

We observe that the actual implementation of the $s$-stage EQUIP method (\ref{Butab1})--(\ref{Xs1}) requires the evaluation of the integrals (\ref{ros12}) and (\ref{gammas12}). Consequently, it is not yet an operative method since, quoting e.g., Dahlquist and Bj\"ork \cite[p.\,521]{DaBj08}, {\em as is well known, even many relatively simple integrals cannot be expressed in finite terms of elementary functions, and thus must be evaluated by numerical methods.} For this purpose, one may use a $k$-point Gauss-Legendre quadrature, and we denote EQUIP$(k,s)$ the corresponding method, whose analysis is now completed, w.r.t. that done in the previous section. Hereafter, we shall assume $k\ge s$.
Consequently, by setting
\begin{equation}\label{cibik}
P_k(\ck_i) = 0, \qquad \bk_i = \int_0^1 \hat\ell_i(x)\dd x, \qquad \hat\ell_i(x) = \prod_{j\ne i} \frac{x-\ck_j}{\ck_i-\ck_j}, \qquad i=1,\dots,k,
\end{equation}
the nodes and weights of the Gauss-Legendre formula of order $2k$, one approximates (\ref{ros12}) and (\ref{gammas12}) by means of
\begin{equation}\label{ros12k}
\hat\rho_j(\sigma_1)  = \sum_{\ell=1}^k \bk_\ell P_j(\ck_\ell) \nabla C(\sigma_1(\ck_\ell h)), \quad j=0,\dots,s-1,\qquad
\hat\rho(\sigma_2) = \sum_{\ell=1}^k \bk_\ell \nabla C(\sigma_2(\ck_\ell)), 
\end{equation} 
and
\begin{equation}\label{gammas12k}
\hat\gamma_j(\sigma_1)  = \sum_{\ell=1}^k \bk_\ell P_j(\ck_\ell) f(\sigma_1(\ck_\ell h)), \quad j=0,\dots,s-1,\qquad
\hat\gamma(\sigma_2) = \sum_{\ell=1}^k \bk_\ell f(\sigma_2(\ck_\ell)), 
\end{equation}
respectively. It is quite straightforward to prove that, under regularity assumptions on $f$ and $C$,
\begin{eqnarray}\nonumber
\hat\rho_j(\sigma_1) & = & \rho_j(\sigma_1) + O(h^{2k-j}),   \\[1mm] \label{erroh}
\hat\gamma_j(\sigma_1) & = & \gamma_j(\sigma_1) + O(h^{2k-j}), \qquad j=0,\dots,s-1,\\[1mm]
\hat\rho(\sigma_2) & = & \bar\rho(\sigma_2) + O(h^{2k}),    \nonumber \\[1mm]
\hat\gamma(\sigma_2) & = & \bar\gamma(\sigma_2) + O(h^{2k}). \nonumber
\end{eqnarray}
As a result, the new approximation is formally still defined by the two polynomials $\sigma_1$ and $\sigma_2$ defined in (\ref{sig1}) and (\ref{sig20})--(\ref{sig2}). However, in such formulae, the parameter 
$\alpha$, defined in (\ref{alfa}) or (\ref{alfaH}), has to be replaced by  an approximated one, say $\hat\aa$, computed by using the quadratures (\ref{cibik})--(\ref{gammas12k}). That is, 
\begin{equation}\label{halfa}
\hat\aa = \frac{ \sum_{j=0}^{s-1} \hat\rho_j(\sigma_1)^\top \gamma_j}{\left( \hat\rho_0(\sigma_1)-\hat\rho(\sigma_2)\right)^\top
\left(\phi_{20}\gamma_0-\phi_{10}\gamma_1\right) +\sum_{j=1}^{s-1} \hat\rho_j(\sigma_1)^\top\left(\phi_{2j}\gamma_0-\phi_{1j}\gamma_1\right)}
\end{equation}
or
\begin{equation}\label{halfaH}
\hat\aa = \frac{ \sum_{j=0}^{s-1} \hat\gamma_j(\sigma_1)^\top J\gamma_j}{\left( \hat\gamma_0(\sigma_1)-\hat\gamma(\sigma_2)\right)^\top J
\left(\phi_{20}\gamma_0-\phi_{10}\gamma_1\right) +\sum_{j=1}^{s-1} \hat\gamma_j(\sigma_1)^\top J\left(\phi_{2j}\gamma_0-\phi_{1j}\gamma_1\right)},
\end{equation}
respectively. The following result can then be proved.
\begin{lem}\label{dalfalem} Let the parameter $\aa$ be defined by either (\ref{halfa}), when solving problem (\ref{fy})--(\ref{Cy}), or (\ref{halfaH}), when solving problem (\ref{H})--(\ref{Hfun}). Then, for all $k\ge s$ one has 
\begin{equation}\label{dalfa}
\hat\aa = \aa + O(h^{2(k-1)}) \qquad \Rightarrow\qquad \hat\aa = O(h^{2(s-1)}),
\end{equation} 
with $\aa$  given by (\ref{alfa}) and (\ref{alfaH}), respectively. In particular, when the problem is Hamiltonian and $k=s$, one has $\hat\aa=0$, so that EQUIP$(s,s)$ reduces to the $s$-stage Gauss-method.
\end{lem}
\proof The first part of the proof follows straightforwardly, by taking into account (\ref{erroh}), from the line integral formulation used to derive (\ref{alfa}) and (\ref{alfaH}), respectively. In particular, one finds that
\begin{equation}\label{halfa1}
\hat\aa = \frac{\N(h) + O(h^{2k})}{\D(h) + O(h^{2k-s+1})} \equiv \aa\cdot \left( 1 + O(h^{2(k-s)})\right),
\end{equation}
with $\N(h) = O(h^{2s})$ and $\D(h)=O(h^2)$ formally given by (\ref{Odih2s}) and  (\ref{D1D2}), in case of problem (\ref{fy})--(\ref{Cy}), or (\ref{Odih2sH}) and (\ref{D1D2H}), in case of problem (\ref{H})--(\ref{Hfun}), respectively.

The fact that $\hat\aa=0$, when the problem is Hamiltonian and $k=s$, follows from (\ref{halfaH}), by considering that matrix $J$ is skew-symmetric and (see (\ref{gammai}) and (\ref{gammas12k})) $\hat\gamma_j(\sigma_1)\equiv \gamma_j$.\QED\medskip

Moreover, the following results hold true, whose straightforward proofs, which strictly follow that of Theorem~\ref{alfaT1}, are omitted for sake of brevity.

\begin{theo}\label{equipksC} With reference to problem (\ref{fy})--(\ref{Cy}) and under regularity assumptions on both $f$ and $C$, and for $k\ge s$, one obtains that the EQUIP$(k,s)$ method defined by the parameter $\hat\aa$ in (\ref{halfa}) satisfies:
$$C(y_1)-C(y_0) = \left\{\begin{array}{cl} 0, & \mbox{ if } C \mbox{ is a polynomial of degree not larger than $2k/s$},\\[2mm]
O(h^{2k+1}), & \mbox{otherwise}.
\end{array}\right.$$
\end{theo}

\begin{theo}\label{equipksH} With reference to problem (\ref{H})--(\ref{Hfun}) and under regularity assumptions on $H$, and for $k\ge s$, one obtains that the EQUIP$(k,s)$ method defined by the parameter $\hat\aa$ in (\ref{halfaH}) satisfies:
$$H(y_1)-H(y_0) = \left\{\begin{array}{cl} 0, & \mbox{ if } H \mbox{ is a polynomial of degree not larger than $2k/s$},\\[2mm]
O(h^{2k+1}), & \mbox{otherwise}.
\end{array}\right.$$
\end{theo}

\begin{rem} As a consequence of Theorems~\ref{equipksC} and \ref{equipksH}, one has that an {\em exact} conservation can be always gained, in case of a polynomial invariant, by taking $k$ large enough. On the other hand, even in the general case a {\em practical} conservation can still be obtained, by taking $k$ large enough so that the quadrature error  is within round-off error level. It is worth mentioning that by considering higher values of $k$ does not increase the complexity of the nonlinear problem to be solved, which amounts to compute the $s$ coefficients of the polynomial $\sigma_1$, and the parameter $\alpha$.
\end{rem}


The previous results allows us to prove that, for all $k\ge s$, the EQUIP$(k,s)$ method still has order $2s$.

\begin{theo}\label{horder} Let $y_1\approx y(h)$ be the approximation provided by the EQUIP$(k,s)$ method (\ref{Butab1})--(\ref{Xs1}), with $\aa=\hat\aa$ given by (\ref{halfa}) (respectively, (\ref{halfaH}) when the problem is Hamiltonian). Then, for all $k\ge s$, and assuming the hypotheses of Theorem~\ref{order} hold true, 
\begin{equation}\label{hord2s}
y_1-y(h) = O(h^{2s+1}), 
\end{equation}
i.e., the method retains the order $2s$.
\end{theo}
\proof The proof strictly follows that done for Theorem~\ref{order}, the only difference being that, starting from (\ref{s1menou}), $\aa$ has to be formally replaced by $\hat\aa$. Nevertheless, because of (\ref{dalfa}), the same arguments of that proof apply also in this case, so that (\ref{hord2s}) follows.\QED\medskip

\begin{rem} It is worth noticing that, by taking into account that $\hat\aa=O(h^{2s-2})$, from (\ref{Ydigamma1}) it follows that the stage order of EQUIP$(k,s)$ methods is $s$, for all $k\ge s$. I.e., the methods also retain the same stage order as that of the underlying $s$-stage Gauss method.
\end{rem}

\section{The nonlinear iteration}\label{fixpoint}

Next, let us duscuss the use of a fixed-point iteration for solving the discrete problem generated by the EQUIP$(k,s)$ method, with $k\ge s$. In principle, this is given by (\ref{iteration1}), with the only difference that we have to consider $\hat\aa$ in place of $\aa$. Therefore, by taking into account (\ref{dalfa})--(\ref{halfa1}), convergence can be proved by using the same arguments in the proof of Theorem~\ref{exists_unique}, provided that $k>s$, which we assume hereafter. There are, however, two improvements over the basic iteration (\ref{exists_unique}), as explained below.

\smallskip
First of all, in order to improve the convergence properties of the iteration, we consider the following {\em Gauss-Seidel type} variant of (\ref{iteration1}),
\begin{equation}
\label{iteration2}
\left\{ 
\begin{array}{l}
\displaystyle \hat\alpha^{\ell+1}=\frac{\hat N(\hat\aa^\ell,\bfgamma^\ell)}{\hat D(\hat\aa^\ell,\bfgamma^\ell)}, \\[.4cm]
\bfgamma^{\ell+1}=\Psi(\hat\aa^{\ell+1},\bfgamma^\ell),  \qquad \ell=0,1,\dots,
\end{array}
\right.
\end{equation}
where we have formally used $\hat N$ and $\hat D$ in place of $N$ and $D$, respectively, in order to take into account of the quadratures (see (\ref{cibik})--(\ref{halfaH})).

\smallskip
Secondly, the quadratures involved in the computation of $\hat\aa$ in place of $\aa$, potentially introduce a $O(h^{2k+1})$ error in the invariant at each step, according to Theorems~\ref{equipksC} and \ref{equipksH}. In order to avoid their accumulation, it suffices to compute $\hat\aa$ so that the variation of the invariant, at the current step, matches the opposite of the error until the previous step. This latter error, which we denote $\Delta C$, can be easily computed by evaluating the invariant. As a result,  following steps similar to those used in the proof of Theorem~\ref{alfaT1}, one obtains that the iteration (\ref{iteration2}) becomes:
\begin{equation}
\label{iteration3}
\left\{ 
\begin{array}{l}
\displaystyle \hat\alpha^{\ell+1}=\frac{\hat N(\hat\aa^\ell,\bfgamma^\ell) + \frac{\Delta C}h}{\hat D(\hat\aa^\ell,\bfgamma^\ell)}, \\[.4cm]
\bfgamma^{\ell+1}=\Psi(\hat\aa^{\ell+1},\bfgamma^\ell),  \qquad \ell=0,1,\dots.
\end{array}
\right.
\end{equation}
In so doing, the errors in the invariant remain uniformly $O(h^{2k+1})$, since their accumulation is prevented, and moreover, the magnitude of $\hat\aa$ remains $O(h^{2(s-1)})$.

\section{Numerical results}\label{numtest}

In this section we consider a few numerical tests concerning the numerical approximation of some conservative problems. In particular, two Hamiltonian problems and two Poisson problems. They are listed below.
\begin{enumerate} 
\item The well-known Kepler problem \cite{HLW06,BI16}, which is Hamiltonian,  i.e. in the form (\ref{H}), with
\begin{equation}\label{kepler}
H(q,p) = \frac1{2} \|p\|_2^2 -\frac{1}{\|q\|_2}, \qquad q,p\in\RR^2.
\end{equation}
It admits also the additional quadratic invariant (i.e., the {\em angular momentum}),
\begin{equation}\label{M} M(q,p) = q^\top J p, \qquad J = \pmatrix{rr} 0 & 1\\-1 & 0\endpmatrix.\end{equation}
Its solution is periodic of period $T=2\pi$, and resulting into an ellipse of eccentricity $\eps$ in the $q$-plane, when starting at 
\begin{equation}\label{kep0}
q(0) = (1-\eps,\,0)^\top, \qquad p(0) = \left(0,\,\sqrt{\frac{1+\eps}{1-\eps}}\right)^\top.
\end{equation}
In particular, we shall consider the value $\eps = 0.5$.

\item The pendulum problem \cite{BI16}, which is also Hamiltonian with
\begin{equation}\label{pendH}
H(q,p) = \frac{1}2p^2-\cos(q).
\end{equation}
Its solution is periodic of period $T\approx  28.57109480185544$, when starting at
\begin{equation}\label{pend0}
q(0) = 0, \qquad p(0) = 1.99999.
\end{equation}

\item The Poisson problem \cite{BCMR12}
\begin{equation}\label{poisson}
\dot y = B(y) \nabla H(y), \qquad  B(y)^\top = -B(y),
\end{equation}
with $ y\in\RR^3$ and
\begin{equation}\label{poisson2}
B(y) = \pmatrix{ccc} 
0 & c_3y_3 &-c_2y_2\\
-c_3y_3 &0& c_1y_1\\
c_2y_2 & -c_1y_1 &0
\endpmatrix, \qquad H(y) = y_1^{12} +\frac{1}2\left[ (y_2-y_3)^2 +(y_1-y_3)^2\right],
\end{equation}
also admitting, as a further invariant besides $H$, the quadratic {\em Casimir} 
\begin{equation}\label{casimir}
C(y) = y^\top \pmatrix{ccc} c_1 \\ &c_2\\ &&c_3\endpmatrix y.
\end{equation}
The solution turns out to be periodic, of period $T\approx0.53102669598427$, when choosing
\begin{equation}\label{cy0}
c_1 = 1, \quad c_2=5,\quad c_3= -4, \qquad\mbox{and}\qquad y(0) = (1,\,1,\,1)^\top.
\end{equation}

\item The {\em Lotka-Volterra} problem, which is in the form (\ref{poisson}) with $y\in\RR^2$ and
\begin{equation}\label{PP}
B(y) = \pmatrix{cc} 
0 &  y_1y_2\\
 -y_1y_2 &0
\endpmatrix, \qquad H(y) = a\log y_1-y_1 +b\log y_2-y_2.
\end{equation}
The solution turns out to be periodic, of period $T\approx7.720315563434113$, when choosing
\begin{equation}\label{aby0}
a = 1, \qquad b =2, \qquad\mbox{and}\qquad y(0) = (0.1,\,0.1)^\top.
\end{equation}
\end{enumerate}
For such problems, we compare the EQUIP$(k,s)$ methods with the corresponding $s$-stage Gauss collocation methods, $s=2,3$,  using a given constant stepsize of the form $h=T/n$, with $T$ the period of the solution and increasing values of $n$. This allows us to easily asses the accuracy of the solution, by checking it at the end of each period. As a measure of the error on the Hamiltonian, by setting $H_i$ its numerical value at the $i$th integration step, we consider the mean square value of the error over $N$ steps,
$$e_H := \sqrt{\frac{1}N\sum_{i=1}^N(H_i-H_0)^2}.$$ A similar measure will be used for the invariants (\ref{M}), $e_M$, and (\ref{casimir}), $e_C$. When needed, by setting $\aa_i$ the computed parameter (\ref{halfa})--(\ref{halfaH}) for the EQUIP method, we set  
$$\bar\aa := \sqrt{\frac{1}N\sum_{i=1}^N\aa_i^2},$$ as a measure of its magnitude. The EQUIP methods are implemented via the iteration (\ref{iteration3}), whereas only the second equation, having fixed $\hat\aa=0$, is used for the Gauss methods.  For the EQUIP$(k,s)$ methods, the value of $k$ is chosen in order to obtain a practical energy-conservation, as soon as the stepsize is decreased, for the considered values of $s$: in the present case, $k=6$ is sufficient for our purposes.

\medskip
We start considering the Kepler problem (\ref{kepler})--(\ref{kep0}). In Table~\ref{kep62} we list the obtained results by solving the problem with the EQUIP(6,2) and the 2-stage Gauss (GAUSS2, hereafter) methods. From such results, one deduces that:
\begin{itemize}
\item both methods are 4-th order, with EQUIP(6,2) one order of magnitude more accurate than GAUSS2, and preserve the quadratic invariant $M$;

\item for EQUP(6,2), $\bar\aa$ turns out to be $O(h^2)$ as expected and, moreover, the error on $H$ soon falls within round-off error level, whereas it decreases with order 4 for GAUSS2;

\item the number of iterations per step (iter/step) is essentially the same for both methods, and it decreases (as is expected) with the stepsize $h$.

\end{itemize}
Similarly, in Table~\ref{kep63} we list the obtained results by solving the problem with the EQUIP(6,3) and the 3-stage Gauss (GAUSS3, hereafter) methods, showing that:
\begin{itemize}
\item both methods are 6-th order, with EQUIP(6,3) one order of magnitude more accurate than GAUSS3, and preserve the quadratic invariant $M$;

\item for EQUP(6,3), $\bar\aa$ turns out to be $O(h^4)$ as expected and, moreover, the error on $H$ soon falls within round-off error level, whereas it decreases with order 6 for GAUSS3;

\item the number of iterations per step (iter/step) is essentially the same for both methods, and it decreases (as is expected) with the stepsize $h$.

\end{itemize}

\begin{table}[th]
\caption{Numerical solution of the Kepler problem (\ref{kepler})--(\ref{kep0}) by using the EQUIP(6,2) and the 2-stage Gauss (GAUSS2) methods with stepsize $h=2\pi/n$ over 10 periods.}
\label{kep62}\smallskip
{\scriptsize
\hspace{-.3cm}\begin{tabular}{|r|rrrrrrr|rrrrrr|}
\hline
                      &\multicolumn{7}{c|}{EQUIP(6,2)} &\multicolumn{6}{c|}{GAUSS2}\\[1mm]
           $n$     & error & rate & $H$-error & $M$-error & $\bar\aa$ & rate & iter/ &  error & rate & $H$-error & rate & $M$-error & iter/ \\ 
                      &          &        &                 &                  &                 &        & step &          &         &                 &        &                  & step \\     
\hline
20  &   1.34e-1  &          --   &1.64e-09 &  3.66e-15 &  1.51e-3       &    --  & 19.6 &  1.55e0    &        --   & 1.95e-3      &      -- & 5.75e-15  & 17.4\\
30  &   2.61e-2  & 4.0  & 6.10-12  & 7.88e-15 &  6.81e-4 &  2.0 &  15.6 &  2.37e-1 &  4.6 &  2.27e-4 &  5.3&   1.15e-14 &  14.2\\
40 &   8.36e-3  & 4.0  & 1.86e-13 &  4.22e-16  & 3.84e-4  & 2.0  & 13.6 &  8.00e-2  & 3.8  & 7.65e-5  & 3.8&  1.24e-15 &  12.8\\
50 &   3.45e-3  & 4.0  & 1.84e-14 &  1.58e-15  &  2.45e-4 & 2.0 &  12.5 &  3.41e-2 &  3.8  & 3.28e-5  & 3.8&   2.39e-15 &  11.7\\
60 &  1.67e-3   & 4.0  & 2.05e-15 &  4.70e-16  &  1.70e-4 &  2.0 &  11.8 &  1.68e-2 &  3.9  & 1.61e-5  & 3.9&   3.16e-15 &  11.3\\
70 &  9.01e-4   & 4.0  & 1.44e-15 &  4.28e-16  &  1.25e-4 &  2.0 &  11.4 &  9.17e-3 &  3.9  & 8.83e-6  & 3.9&   4.10e-15 &  10.6\\
80 &  5.29e-4   & 4.0  & 1.11e-15 &  2.46e-15  &  9.58e-5 &  2.0  & 10.8 &  5.41e-3 &  3.9  & 5.22e-6  & 3.9&   1.36e-15 &  10.3\\
90 &  3.31e-4   & 4.0  & 2.44e-15  & 8.11e-16  &  7.57e-5 &  2.0 &  10.5 &  3.40e-3 &  4.0  & 3.27e-6  & 4.0&   3.33e-15 &  10.0\\
100 & 2.18e-4   & 4.0 &  9.77e-16 &  2.98e-15  & 6.13e-5 &  2.0 &  10.2 &  2.24e-3 &  4.0  & 2.16e-6  & 4.0&  6.65e-15 &  9.7\\
\hline
\end{tabular}}
\caption{Numerical solution of the Kepler problem (\ref{kepler})--(\ref{kep0}) by using the EQUIP(6,3) and the 3-stage Gauss (GAUSS3) methods with stepsize $h=2\pi/n$ over 10 periods.}
\label{kep63}\smallskip
{\scriptsize
\hspace{-.3cm}\begin{tabular}{|r|rrrrrrr|rrrrrr|}
\hline
                      &\multicolumn{7}{c|}{EQUIP(6,3)} &\multicolumn{6}{c|}{GAUSS3}\\[1mm]
           $n$     & error & rate & $H$-error & $M$-error & $\bar\aa$ & rate & iter/ &  error & rate & $H$-error & rate & $M$-error & iter/ \\          
                      &          &        &                 &                  &                 &        & step &          &         &                 &        &                  & step \\     
\hline
20  &   2.67e-3   &        --  & 1.15e-09 &  5.44e-15  & 4.62e-5  &    -- &  15.3  & 5.16e-2   & -- &   6.72e-5    & -- &  1.55e-15  & 15.4\\
30  &   3.11e-4 &  5.3 &   1.68e-11  & 4.43e-15  & 1.17e-5 &  3.4 &  13.1 &  7.41e-3 &  4.8 &  8.44-6 &  5.1 & 2.17e-15 &  13.1 \\
40 &   5.63e-5 &  6.0  & 4.61e-13  & 2.00e-15  & 3.81e-6  & 3.9 &  11.9 &  1.22e-3 &  6.3 &  1.38e-6 &  6.3  & 4.08e-15  & 11.9\\
50 &   1.47e-5  & 6.0 &  2.38e-14  & 9.81e-16  & 1.55e-6 &  4.0 &  11.3 &  3.09e-4 &  6.2 &  3.48e-7 &  6.2  & 8.59e-16 &  11.3\\
60 &  4.94e-6  & 6.0  & 2.01e-15  & 7.91e-16  & 7.47e-7 &  4.0  & 10.5 &  1.02e-4 &  6.1 &  1.15e-7 &  6.1  & 9.51e-16 &  10.5\\
70 &  1.96e-6  & 6.0   & 7.58e-16 &  5.67e-16  & 4.02e-7 &  4.0 &  10.1 &  4.01e-5 &  6.0 &  4.51e-8 &  6.1  & 6.47e-16  & 10.1\\
80 &  8.78e-7 &  6.0  &  3.53e-16  & 8.32e-16  & 2.35e-7 &  4.0 &  9.7 &  1.79e-5 &  6.0 &  2.01e-8 &  6.0 & 9.04e-16 &  9.7\\
90 &  4.33e-7  & 6.0  &  3.92e-16  & 1.13e-15  & 1.47e-7 &  4.0 &  9.3 &  8.82e-6 &  6.0 &  9.90e-9 &  6.0 &  2.22e-15 &  9.3\\
100 & 2.30e-7  & 6.0  & 2.46e-16 &  4.30e-16  & 9.62e-8 &  4.0 &  9.1 &  4.68e-6 &  6.0 &  5.25e-9 &  6.0 &  6.18e-16 &  9.1\\
\hline
\end{tabular}}
\end{table}
\begin{table}[th]
\caption{Numerical solution of the pendulum problem (\ref{pendH})--(\ref{pend0}) by using the EQUIP(6,2) and the 2-stage Gauss (GAUSS2) methods with stepsize $h=T/n$ over 10 periods.}
\label{pen62}\smallskip
{\scriptsize 
\centerline{\begin{tabular}{|r|rrrrrr|rrrrr|}
\hline
                      &\multicolumn{6}{c|}{EQUIP(6,2)} &\multicolumn{5}{c|}{GAUSS2}\\[1mm]
           $n$     & error & rate & $H$-error & $\bar\aa$ & rate & iter/  &  error & rate & $H$-error & rate & iter/ \\ 
                      &          &        &                 &                  &        & step &          &         &                 &        & step \\     
\hline
50   &   4.54e-1 &    --  &  4.71e-06 &  4.55e-3 &  --   &   40.6  & 1.48e2  & --      & 2.61e-4  &  --   &   17.5\\ 
60   &   2.23e-1 &  3.9 &  1.38e-11  & 3.22e-3  & 1.9 &   34.4  & 1.35e2  & 0.5   &  1.05e-4 &  5.0 &   16.1\\
70   &   1.22e-1 &  3.9 &  7.47e-12 &  2.40e-3  & 1.9 &   32.5  & 1.23e2  & 0.6   & 5.16e-5  & 4.6  &   15.0\\
80   &   7.22e-2 &  3.9 &  6.78e-09 &  1.84e-3  & 2.0 &   30.4  & 2.86e0  & 28 &  2.72e-5 &  4.8 &   14.0\\
90   &  4.54e-2  & 3.9  & 9.33e-13  & 1.50e-3  & 1.8  &   27.8  & 3.63e0  & -2.0  &  1.72e-5 &  3.9 &   13.4\\
100 &  3.01e-2  & 3.9  & 4.73e-13  & 1.22e-3  & 2.0  &   25.2  & 3.72e0  & -0.2 &  1.12e-5 &  4.1 &   12.9\\
110 &  2.09e-2  & 3.8  & 2.46e-13  & 1.01e-3  & 1.9  &   23.5  & 3.70e0  &  0.0  & 7.53e-6  &  4.2 &    12.4\\
120 &  1.52e-2  & 3.7  & 1.03e-13  & 8.71e-4  & 1.8  &   23.3  & 3.58e0  &  0.4  & 5.25e-6  & 4.1  &   12.0\\
130 &  1.11e-2  & 3.9  & 5.51e-14  & 7.48e-4  & 1.9  &    22.0 &  3.29e0 &  1.1  & 3.78e-6  & 4.1  &    11.6\\
140 &  8.36e-3  & 3.9  & 3.34e-14  & 6.50e-4  & 1.9  &    20.5 &  2.86e0 &  1.9  & 2.79e-6  & 4.1  &   11.4\\
150 &  6.31e-3  & 4.1  & 2.49e-14  & 5.65e-4  & 2.0  &    20.5 &  2.37e0 &   2.7 &  2.10e-6 &  4.1 &   11.2\\
\hline
\end{tabular}}}
\caption{Numerical solution of the pendulum problem (\ref{pendH})--(\ref{pend0}) by using the EQUIP(6,3) and the 3-stage Gauss (GAUSS3) methods with stepsize $h=T/n$ over 10 periods.}
\label{pen63}\smallskip
{\scriptsize
\centerline{\begin{tabular}{|r|rrrrrr|rrrrr|}
\hline
                      &\multicolumn{6}{c|}{EQUIP(6,3)} &\multicolumn{5}{c|}{GAUSS3}\\[1mm]
           $n$     & error & rate & $H$-error & $\bar\aa$ & rate & iter/  &  error & rate & $H$-error & rate & iter/ \\ 
                      &          &        &                 &                  &        & step &          &         &                 &        & step \\     
\hline
50   &   7.67e-3 &  --   &  9.20e-12 &  2.0972e-4 &  --    &  16.4 &  3.72e0  &     -- &   1.01e-5&  --   &   15.8\\ 
60   &   1.80e-3 & 7.9 &  6.04e-12 &  6.6569e-5 &  6.3  &  14.7 &  3.36e0  &  0.6 &  2.86e-6 &  6.9 &   14.6\\
70   &   5.69e-4 & 7.5 &  3.75e-12 &  2.7582e-5 &  5.7  &  13.9 &  1.97e0  &  3.5 &  1.07e-6 &  6.4 &   13.8\\
80   &   2.16e-4 & 7.2 &  2.30e-12 &  1.3510e-5 &  5.3  &  13.1 &  9.29e-1 &  5.6 &  4.65e-7 &  6.2 &   13.2\\
90   &  1.28e-4 &  4.5 &  1.13e-12 &  9.2616e-6 &  3.2  &  12.6 &  4.56e-1 &  6.0 &  2.24e-7 &  6.2 &   12.6\\
100 &  6.19e-5 &  6.9 &  8.93e-13 &  5.6190e-6 &  4.7  &  12.1 &  2.40e-1 &  6.1 &  1.18e-7 &  6.1 &   12.3\\
110 &  3.15e-5 &  7.1 &  5.72e-13 &  3.5084e-6 &  4.9  &  11.8 &  1.34e-1 &  6.1 &  6.57e-8 &  6.1 &    11.8\\
120 &  1.72e-5 &  7.0 &  3.74e-13 &  2.3050e-6 &  4.8  &  11.4 &  7.92e-2 &  6.1 &  3.87e-8 &  6.1 &   11.5\\
130 &  9.83e-6 &  7.0 &  2.50e-13 &  1.5774e-6 &  4.7  &  11.1 &  4.88e-2 &  6.1 &  2.38e-8 &  6.1 &    11.2\\
140 &  5.88e-6 &  6.9 &  1.70e-13 &  1.1162e-6 &  4.7  &  10.9 &  3.11e-2 &  6.1 &  1.52e-8 &  6.1 &    10.9\\
150 &  3.65e-6 &  6.9 &  1.22e-13 &  8.1196e-7 &  4.6  &  10.6 &  2.05e-2 & 6.0  &  1.00e-8 & 6.1  &   10.6\\
\hline
\end{tabular}}}
\end{table}

\medskip
Next, we briefly study the numerical solution of the pendulum problem (\ref{pendH})--(\ref{pend0}) by means of EQUIP$(k,s)$  and the corresponding Gauss  method. 
In particular:
\begin{itemize}
\item in Table~\ref{pen62} we list the obtained results for the EQUIP(6,2) and GAUSS2 methods, showing that, even though GAUSS2 is able to conserve the Hamiltonian with the prescribed order, nonetheless, the EQUIP(6,2) method is much more accurate in conserving the Hamiltonian, and this reflects into a much more accurate numerical solution;

\item in Table~\ref{pen63} we list the obtained results for the EQUIP(6,3) and GAUSS3 methods, showing that also in this case, the energy-conserving EQUIP(6,3) is much more reliable than the GAUSS3 method.

\end{itemize}

\begin{figure}[p]
\centerline{\includegraphics[width=7.5cm,height=6cm]{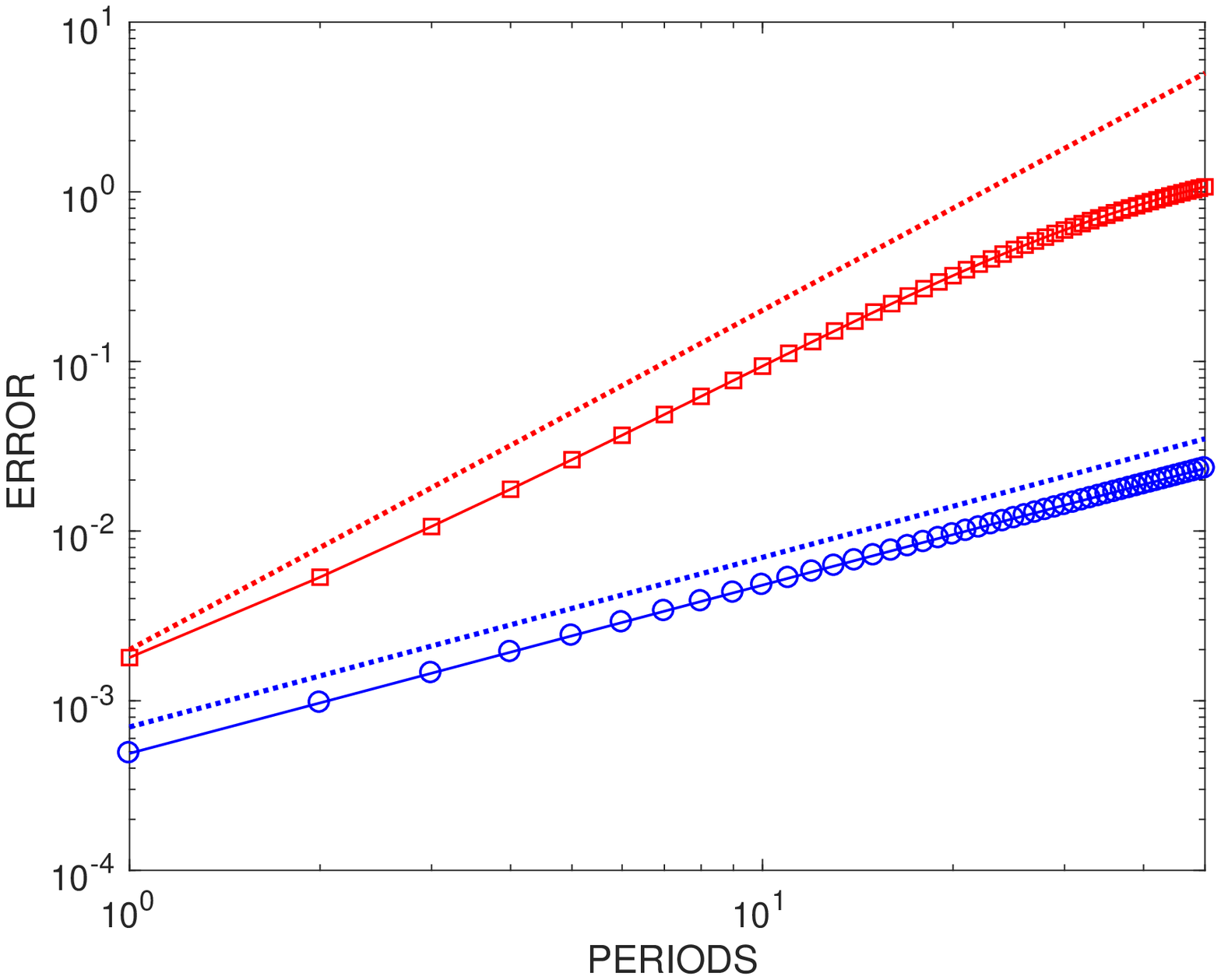} \qquad \includegraphics[width=7.5cm,height=6cm]{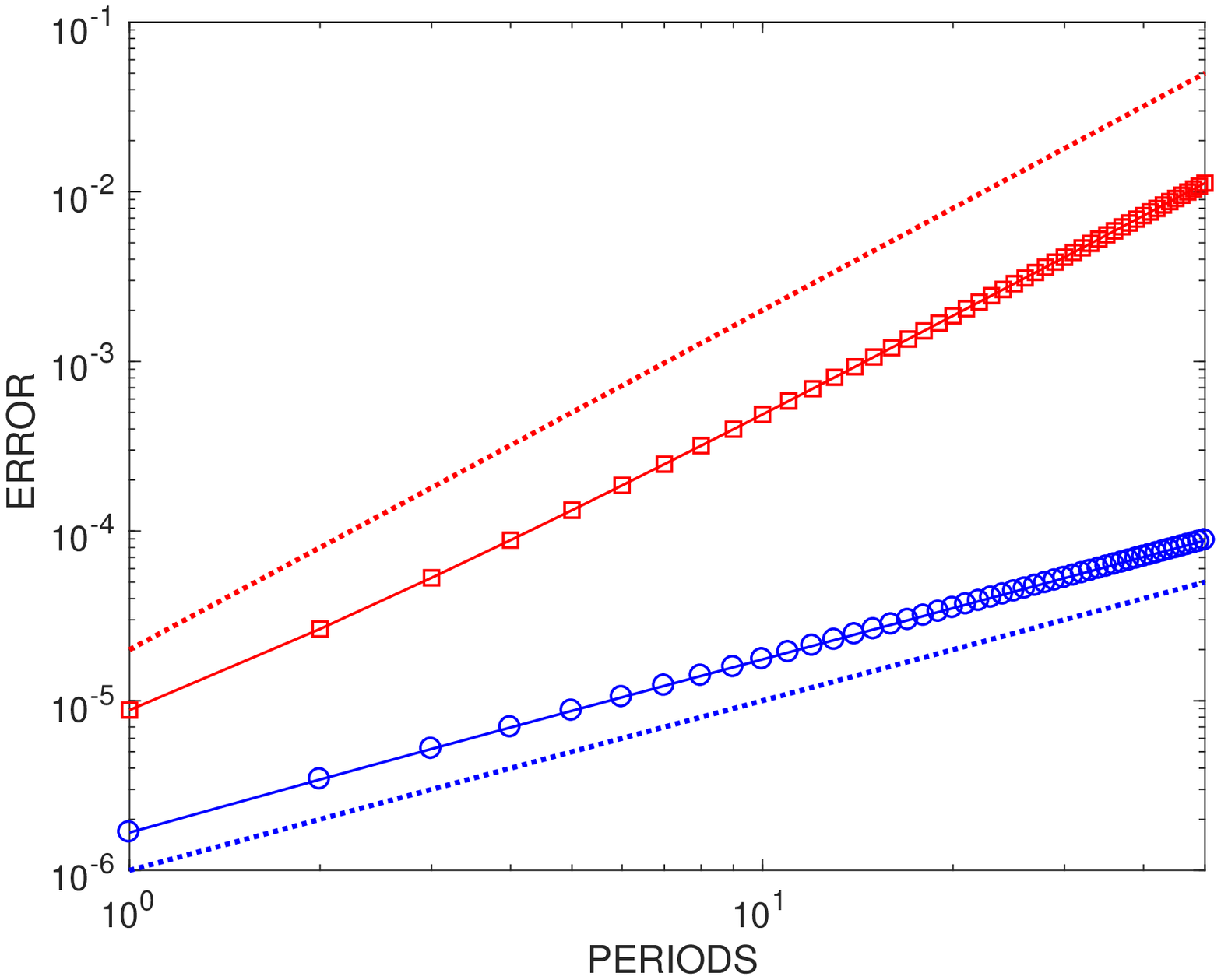}}
\caption{Poisson problem (\ref{poisson})--(\ref{cy0}). Left-plot: solution error along 50 periods for EQUIP(6,2) (circles and solid line)  and GAUSS2 (squares and solid line).   Right-plot: solution error along 50 periods for EQUIP(6,3) (circles and solid line)  and GAUSS2 (squares and solid line). In both plots, the linear and quadratic grows are reported in dotted-line.}
\label{poissonfig}
\bigskip
\centerline{\includegraphics[width=7.5cm,height=6cm]{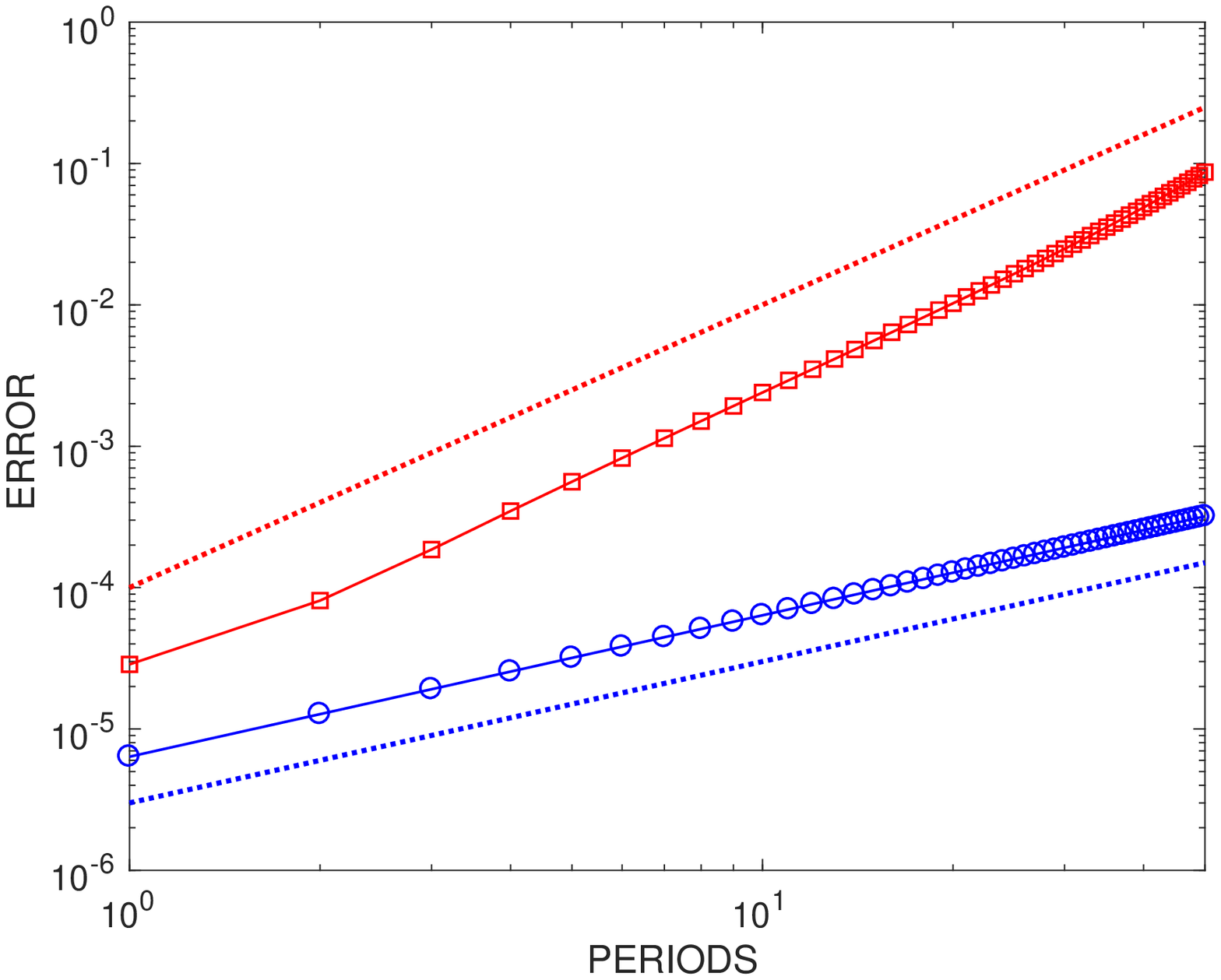} \qquad \includegraphics[width=7.5cm,height=6cm]{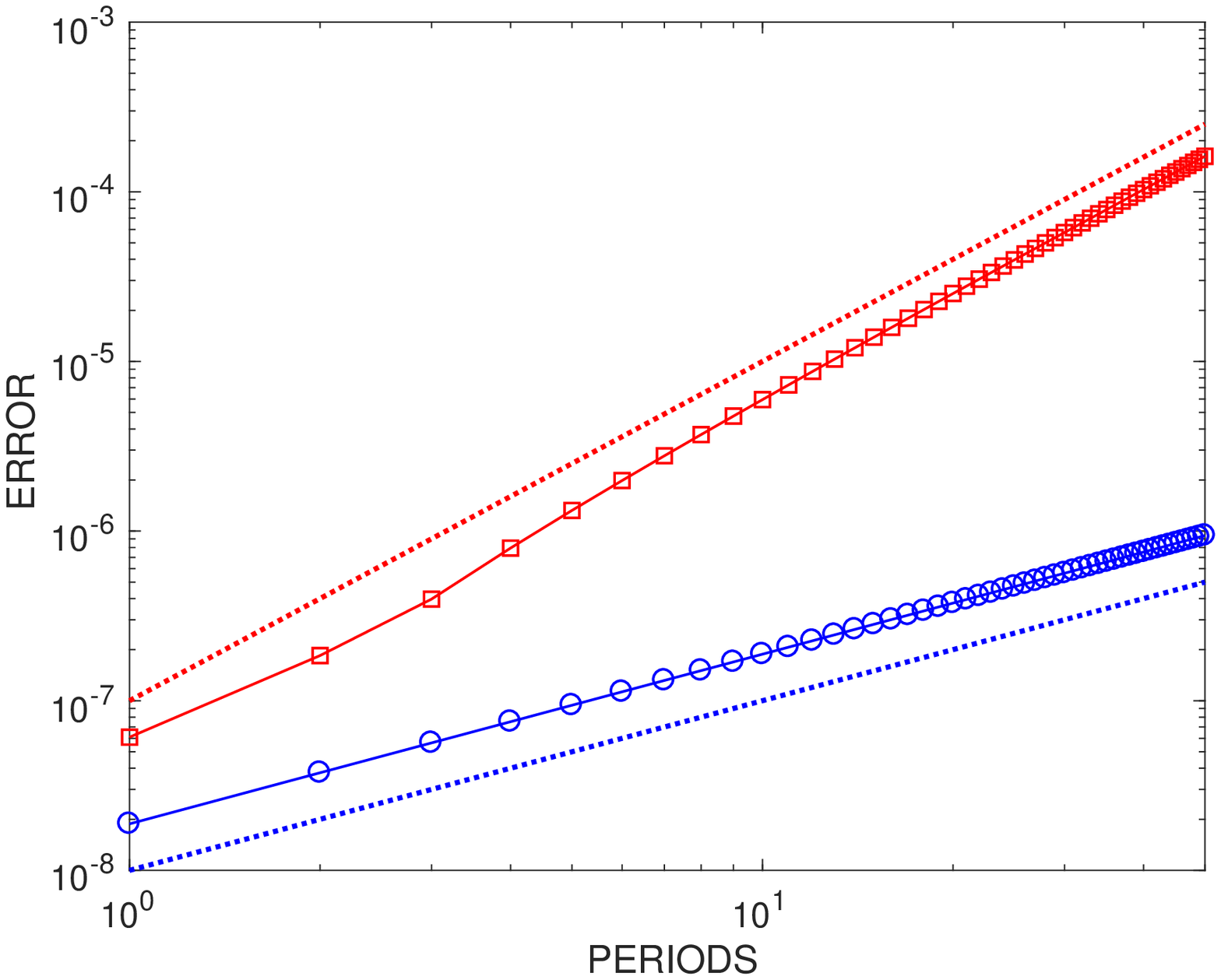}}
\caption{Lotka-Volterra problem (\ref{PP})--(\ref{aby0}).  Left-plot: solution error along 50 periods for EQUIP(6,2) (circles and solid line)  and GAUSS2 (squares and solid line).   Right-plot: solution error along 50 periods for EQUIP(6,3) (circles and solid line)  and GAUSS2 (squares and solid line). In both plots, the linear and quadratic grows are reported in dotted-line.}
\label{ppfig}
\end{figure}

\medskip
We now consider the numerical solution of the Poisson problem (\ref{poisson})--(\ref{cy0}) by means of the EQUIP(6,2)  and GAUSS2 methods, and the EQUIP(6,3)  and  GAUSS3 methods, over 50 periods, by using a stepsize $h=T/100$. In Figure~\ref{poissonfig}  are the computed results for the 4th-order (left-plot) and 6th-order (right-plot) methods, all of them exactly conserving the quadratic Casimir (\ref{casimir}). As one may see, the EQUIP methods exhibit a linear error growth, whereas the Gauss methods show a quadratic error growth.

\medskip
A similar result is obtained when numerically solving the Lotka-Volterra problem (\ref{PP})--(\ref{aby0}) by means of the EQUIP(6,2)  and GAUSS2 methods, and the EQUIP(6,3)  and  GAUSS3 methods, over 50 periods, by using a stepsize $h=T/100$. In Figure~\ref{ppfig}  are the computed results for the 4th-order (left-plot) and 6th-order (right-plot) methods. As one may see, also in this case the EQUIP methods exhibit a linear error growth, whereas the Gauss methods show a quadratic error growth.

\section{Conclusions}\label{fine}

In this paper we reported a novel analysis for EQUIP($k,s$) methods, a symplectic variant of the $s$-stage Gauss collocation methods with energy-preserving properties. Such methods, which retain the order $2s$ of the underlying Gauss method, were originally devised for Hamiltonian problems and are here extended to cope with general conservative problems. The practical implementation of the methods has been another major focus of the paper, ranging from the discretization of the involved integrals, by means of suitable high-order Gaussian rules, to the solution of the generated discrete problem. In particular, the fixed-point iteration used for the analysis of the methods has been improved, in order to maintain a uniform error level in the invariant. A few numerical tests, on a couple of Hamiltonian and Poisson problems, confirm the theoretical findings and show that such methods can be far more reliable than the original Gauss formulae, due to their better conservation properties. 

As a future direction of investigation we mention that, upon devising a Newton-type iteration more efficient than (\ref{iteration2}), EQUIP($k,s$) methods could be conveniently applied also for numerically solving Hamiltonian PDEs, after a proper space semi-discretization. In such a case, in fact, they will conserve possible quadratic invariants of the solution, besides the (semi)-discrete Hamiltonian (see, e.g., \cite{BBGCI17}).

\end{document}